\documentclass{elsarticle}
\usepackage{graphics,xcolor}
\usepackage[normalem]{ulem}
\usepackage{amssymb}
\usepackage{amsmath}
\usepackage{amsfonts}
\usepackage{mathrsfs}%
\usepackage{epsfig,color}
\usepackage{multirow}%
\usepackage{tikz}
\usepackage[mathscr]{euscript}
\usepackage{epstopdf}
\usepackage{ntheorem}
\usepackage{caption}
\usepackage{wrapfig}
\usepackage{array,multirow}
\usepackage{cases}
\usepackage{algpseudocode}
\usepackage[ruled,vlined]{algorithm2e}
\usepackage{tabulary}
\usepackage{breqn}
\usepackage{lmodern}
\usepackage{relsize}
\usepackage{subfig}
\usepackage{enumitem}
\usepackage{physics}
\usepackage{bm}
\usepackage{makecell}
\usepackage{hyperref}

\newcommand{\eps}{\varepsilon}

\newtheorem{theorem}{Theorem}[section]

\newtheorem{remark}[theorem]{Remark}

\newtheorem{example}[theorem]{Example}

\catcode`@=11 \@addtoreset{equation}{section}
\renewcommand\theequation{\thesection.\@arabic\c@equation}
\catcode`@=12

\makeatother

\journal{Elsevier}

\begin{document}
	
\begin{frontmatter}
\title{A Parameter-Driven Physics-Informed Neural Network Framework for Solving Two-Parameter Singular Perturbation Problems Involving Boundary Layers}

\author[address1]{Pradanya Boro}
\ead{pradanyaboro@kgpian.iitkgp.ac.in}
\author[address2]{Aayushman Raina}
\ead{aayushman.raina@iitg.ac.in}
\author[address2]{Srinivasan Natesan \corref{correspondingauthor}}
\ead{natesan@iitg.ac.in}
\cortext[correspondingauthor]{Corresponding author}	
\address[address1]{Department of Mathematics, Indian Institute of Technology Kharagpur, Kharagpur, West Bengal, 721302 - {\sc India}}
\address[address2]{Department of Mathematics, Indian Institute of Technology Guwahati, Guwahati, 781 039 - {\sc India}}

\begin{abstract}
In this article, our goal is to solve two-parameter singular perturbation problems (SPPs) in one- and two-dimensions using an adapted Physics-Informed Neural Networks (PINNs) approach. Such problems are of major  importance in engineering and sciences as it appears in control theory, fluid and gas dynamics, financial modelling and so on. Solutions of such problems exhibit boundary and/or interior layers, which make them difficult to handle. It has been validated in the literature that standard PINNs have low accuracy and can't handle such problems efficiently. Recently Cao et. al \cite{cao2023physics} proposed a new parameter asymptotic PINNs (PA-PINNs) to solve one-parameter singularly perturbed convection-dominated problems. It was observed that PA-PINNs works better than standard PINNs and gPINNs in terms of accuracy, convergence and stability. In this article, for the first time robustness of PA-PINNs will be validated for solving two-parameter SPPs.
\end{abstract}

\end{frontmatter}

\goodbreak\noindent
\section{Introduction}\label{secint}
Parameter dependent differential equations arise in modeling of several physical problems like fluid and gas dynamics, semiconductor devices, elasticity and many more. When the parameter value becomes very small then it leads to singular perturbation problem (SPP). Now due to the influence of small parameters the solution of SPPs exhibits boundary or interior layers. There are several research works carried out to solve SPPs for {\em e.g.},  using finite difference methods \cite{singh2020parameter}, standard finite element method \cite{zhang2021high}, Weak Galerkin method \cite{toprakseven2023weak, raina2024wg}. However, because of the boundary layer(s), standard numerical methods struggle to accurately capture the solution profile, are inefficient, and produce non-physical oscillations. For more details the interested readers may refer to Miller et al. \cite{miller_book_1996} and references therein.
\textcolor{black}{
There is a wealth of research available in the literature on one-parameter problems and one can refer to \cite{hemker2000varepsilon,mohapatra2010parameter, majumdar2017alternating, raina2024fourth} for additional details.}

\textcolor{black}{For two-parameter problems, the solution structure is more intricate than in one-parameter problems, since the solution profile depends on the ratio of the two perturbation parameters (see Section \ref{model_problem_two_para}) \cite{o2006numerical}. This dependency makes numerical solutions more difficult to obtain, compared to one-parameter problems.} Due to the importance of two-parameter SPPs in engineering and scientific fields like lubrication theory \cite{diprima1968asymptotic}, DC motor analysis \cite{chen1974asymptotic}, lot of numerical methods have been developed to solve such problems. \textcolor{black}{Cheng \cite{cheng2021local} studied the local discontinuous Galerkin (LDG) method for a 1D two-parameter SPPs.  Das and Natesan in \cite{avijit2022convergence} considered time-dependent two-parameter SPPs and implemented Streamline diffusion finite element method (SDFEM) in space and Crank-Nicolson scheme in time. Brdar et al. \cite{brdar2016singularly} implemented the standard Galerkin method for 2D two-parameter SPPs on a Bakhvalov-type mesh. In \cite{o2011parameter}, O'Riordan and Pickett proved the parameter-uniform convergence of their numerical scheme applied on a class of two-parameter SPPs in 2D over a Shishkin mesh, which was based on an upwind finite difference operator. Barman et al. developed an ADI-type scheme for 2D parabolic two-parameter SPPs \cite{barman2023alternating, barman2024parameter}. Some more works on solving two-parameter SPPs can be found in \cite{linss2004analysis,singh2020study,raina2025novel}.}

It is very important to have knowledge about the width and location of the boundary layer(s) before the implementation of a numerical technique. Such problems in higher dimensions make the computation very difficult and costly. \textcolor{black}{Contrary to such standard numerical methods, deep neural networks are mesh-less and so does not require to generate a mesh according to the boundary layer location or width, can efficiently manage higher-dimensional problems and exhibit robust learning capabilities. However, there is still a gap concerning the study and application of deep learning methods to solve such problems.}

Deep neural networks exhibit the universal approximation capability, which allows them to estimate any continuous or measurable finite-dimensional function with arbitrary accuracy \cite{hornik1991approximation}. Consequently, they have become a popular tool for solving partial differential equations (PDEs) and Physics Informed Neural Network (PINN) is one of them \cite{raissi2019physics}. The idea is to minimize the loss function which is obtained from the PDE and the associated boundary/initial conditions. It is a mesh-less method and so can tackle any irregular domain easily. Due to its function approximation abilities and simple formulation, PINNs have been employed to solve many different kinds of PDEs which includes stochastic PDEs \cite{zhang2020learning}, integro-differential equations \cite{pang2019fpinns}, fractional equations \cite{pang2019fpinns} to name a few. Problems which have non-smooth solutions can also be tackled by a variant of PINNs which uses the variational formulation of the given problem and is known as $hp$-VPINN \cite{kharazmi2021hp}. Traditional PINNs have many limitations too. It fails to capture the solution profile properly in case of thin boundary layers and steep gradients. It suffers from the spectral bias, {\em i.e.}, the network typically learns broader features before finer details, potentially leading to a skewed representation of the solutions and hence convergence issue arise. In recent years, various modifications of the original PINNs approach have been introduced to address their limitations. These include Conservative PINNs (cPINNs) \cite{jagtap2020conservative}, Finite basis PINNs (FB-PINNs) \cite{moseley2023finite}, Self adaptive PINNs (SA-PINNs) \cite{mcclenny2023self}, Gradient Enhanced PINNs (gPINNs) \cite{yu2022gradient}, and Extended PINNs (XPINNs) \cite{jagtap2021extended}, among others, but none of them addressed thin boundary layers. \textcolor{black}{In \cite{shukla2021parallel}, a general approach for parallelizing PINNs with domain decomposition for flow problems is discussed. Furthermore, additional studies on PINNs have emphasized advancements in neural network architecture and training techniques, {\em e.g.}, adaptive activation function \cite{jagtap2020adaptive} and Convolutional neural network (CNN) type neural architecture \cite{gao2021phygeonet}.}

Recently, Ben et al. \cite{moseley2023finite} proposed finite basis PINNs (FB-PINNs) to address the spectral bias. They showed that FB-PINNs works well for multi-scale problems and in terms of accuracy it surpasses the standard PINNs. Cao et al. \cite{cao2023physics} studied a new parameter asymptotic PINNs (PA-PINNs) to solve one parameter singularly perturbed convection-dominated problems. The idea is to first approximate the solution in the region away from the boundary layer, by inputting large value of the perturbation parameter, optimizing the network for these values and then using these network parameters to train the neural network for small perturbation parameter values and hence approximating the boundary layer solution.

To the best of our knowledge till date PA-PINNs have not been studied for two-parameter SPPs and hence, to move forward in this direction, \textcolor{black}{the aim of this article is to evaluate the performance of PA-PINNs for one- and two-dimensional two-parameter SPPs, compare the results with standard numerical techniques, and demonstrate the advantages of PA-PINNs over these methods}. \textcolor{black}{Novelty of this work lies in generalizing the applicability of the PA-PINNs method by addressing more complex SPPs \cite{o1967two}. The proposed algorithm effectively handles various solution regimes without altering the network architecture, presents a viable and smart approach for solving two-parameter SPPs, as it does not require prior knowledge of the boundary layer location and width. By progressively refining the parameters and iteratively training the PINN, the network's characterization capability is enhanced, which in turn improves the accuracy and convergence of the PINN for approximating problems with large gradients.} 
We have organised the rest of this article as follows: In Section \ref{model_problem_two_para}, we give the model problem setup in one- and two-dimensions. In Section \ref{modi_pinn}, the description of fully connected neural network along with standard PINNs is given. Further, the detailed description of PA-PINNs methodology is also given in Section \ref{modi_pinn}. In Section \ref{num_ex}, numerical results are provided which tells about the performance of PA-PINNs for considered examples. Lastly, in Section \ref{conclusion} we have concluded our work.

\section{Model Problems}\label{model_problem_two_para}

\subsection{One-dimensional two-parameter singular perturbation problem}

\subsubsection{Time-independent case:}

Consider the following two-point boundary-value problem (BVP):
\begin{equation}\label{(ODE)}
\left\{
\begin{array}{ll}
{L}\mathfrak{v}(x)\equiv -\varepsilon_1 \mathfrak{v}''(x)+ \varepsilon_2  \mathfrak{d}(x)\mathfrak{v}'(x)  + \mathfrak{r}(x)\mathfrak{v}(x) = h(x),
\quad x\, \in I = (0,1), \\ [6pt]
\mathfrak{v}(0)  = \mathfrak{b}_{1}, \quad \mathfrak{v}(1) = \mathfrak{b}_{2},
\end{array}\right.
\end{equation}
where $\varepsilon_1 >0 $ and $\varepsilon_2$ are small parameters, $\mathfrak{d}$, $\mathfrak{r}$ and $h$ are sufficiently smooth functions, such that $\mathfrak{d}(x) \geq \gamma > 0$ on $\overline{I} = [0,1]$. The given problem has a unique solution and at both the ends of the domain {\em i.e.}, at $x=0$ and $x=1$ the solution has exponential type of boundary layers \cite{malley_book_1974}.

In order to see the boundary layer behaviour, let $\rho_0$ and $\rho_1$ be the two solutions of the characteristic equation corresponding to \eqref{(ODE)}:
\begin{equation*}
    -\varepsilon_1 \rho^2(x)+ \varepsilon_2  \mathfrak{d}(x)\rho(x)  + \mathfrak{r}(x) = 0.
\end{equation*}
$\rho_1(x) > 0$ and $\rho_0(x) < 0$ describe the layer phenomena near $x=1$ and $x=0$, respectively. Fixing
\begin{equation*}
    \textcolor{black}{\mu_L} = -\max_{x\in[0,1]} \rho_{0}, \quad  \textcolor{black}{\mu_R} = \min_{x\in[0,1]} \rho_{1},
\end{equation*}
or, equivalently,
\begin{equation}\label{mu}
    \textcolor{black}{\mu_{L,R}} = \min_{x\in [0,1]}\dfrac{\mp \varepsilon_2 \mathfrak{d}(x) + \sqrt{\varepsilon_2^2\mathfrak{d}(x)^2 + 4\varepsilon_1\mathfrak{r}(x)}}{2\varepsilon_1}.
\end{equation}
The three regimes based on the relationship between $\varepsilon_1$ and $\varepsilon_2$ are given as follows:
\begin{itemize}
\item If $\varepsilon_1 \ll \varepsilon_2 = 1$, then $\mu_L = \mathcal{O}(1)$ and $\mu_R = \mathcal{O}(1/\varepsilon_1)$. In this case \eqref{(ODE)} is of convection-diffusion type.
\item If $\varepsilon_1 \ll \varepsilon_2^2 \ll 1$, then $\mu_L = \mathcal{O}(1/\varepsilon_2)$ and $\mu_R = \mathcal{O}(\varepsilon_2/\varepsilon_1)$ and is a case when \eqref{(ODE)} behaves as diffusion-convection-reaction type.
\item If $\varepsilon_2^2 \ll \varepsilon_1 \ll 1$, then $\mu_L = \mathcal{O}(1/\sqrt{\varepsilon_1})$ and $\mu_R = \mathcal{O}(1/\sqrt{\varepsilon_1})$ and \eqref{(ODE)} is of reaction-diffusion type.
\end{itemize}

\subsubsection{Time-dependent case:}
Here, we consider the following parabolic problem:
	\begin{equation}\label{ch3_model}
		\left\{
		\begin{array}{lll}
			\eps_1 \dfrac{\partial^2 \mathfrak{v}}{\partial x^2} + \eps_2 \mathfrak{d}(x,t) \dfrac{\partial \mathfrak{v}}{\partial x} - \mathfrak{r}(x,t) \mathfrak{v}(x,t) - \dfrac{\partial \mathfrak{v}}{\partial t} = h(x,t), \ \ (x,t)\in \mathfrak{Q},\\ [10pt]
			\mathfrak{v}(0,t) = \mathfrak{v}_0(t), \ \ \mathfrak{v}(1,t) = \mathfrak{v}_1(t), \quad 0\leq t\leq \mathscr{T},\\ [6pt]
			\mathfrak{v}(x,0) = \phi(x), \ \ x\in \overline{I},
		\end{array} \right.
	\end{equation}	
	where $0 < \varepsilon_1 \ll 1$ and $0 < \varepsilon_2 \ll 1$ are two small parameters, $I = (0,1)$ and $\mathfrak{Q} = I\times (0, \mathscr{T}]$. We assume the boundary and initial data exhibit adequate smoothness, along with compatibility at the corners, to ensure the existence of a unique solution. As $\varepsilon_1$ and $\varepsilon_2$ approach zero, the solution develops layers at the two end points of the domain $I$.

\subsection{Two-dimensional two-parameter singular perturbation problem}

\subsubsection{Elliptic boundary-value problem}

Consider the following 2D elliptic boundary-value problem:
\begin{equation}\label{elliptic_problem}
\left\{
\begin{array}{lll}
\mbox{{L}}_{\varepsilon_1,\varepsilon_2} \mathfrak{v} = h(x,y), &(x,y)\in\varOmega = (0,1)\times(0,1),\,\,\\ [6pt]
\mathfrak{v}(x,y) = 0, &(x,y) \in \partial \varOmega = \overline{\varOmega} \backslash \varOmega,
\end{array} \right.
\end{equation}
where  $0 < \varepsilon_1, \varepsilon_2 \ll 1$ are two small perturbation parameters and the operator $\mbox{{L}}_{\varepsilon_1,\varepsilon_2}$ can be defined as
\begin{equation}\label{ch1_spat_oprtr_2}
	\mbox{{L}}_{\varepsilon_1,\varepsilon_2} \mathfrak{v} := \varepsilon_1 \Delta \mathfrak{v} + \varepsilon_2 \mathfrak{d}(x,y) \cdot \nabla \mathfrak{v} - \mathfrak{r}(x,y) \mathfrak{v}
\end{equation}
or
\begin{equation}\label{ch1_spat_oprtr}
	\mbox{{L}}_{\varepsilon_1,\varepsilon_2} \mathfrak{v} := - \varepsilon_1 \Delta \mathfrak{v} + \varepsilon_2 \mathfrak{d}_1(x,y)\mathfrak{v}_x + \mathfrak{r}(x,y) \mathfrak{v}.
\end{equation}
We define the convection coefficient $\mathfrak{d}$ as $\mathfrak{d}(x,y) = \big( \mathfrak{d}_1(x,y), \,\ \mathfrak{d}_2(x,y) \big)$ such that $\mathfrak{d}_1 \geq \beta_1 > 0$, $\mathfrak{d}_2 \geq \beta_2 > 0$ and $\mathfrak{r}(x,y) \geq r_0 > 0 $. We assume that $\mathfrak{d}, \mathfrak{r}$ and $h$ are smooth enough and $h$ meets the following compatibility property:
\begin{align}\label{ch1_compatible}
h(0,0) = h(0,1) = h(1,0) = h(1,1) = 0.
\end{align}
For the case when $\eps_2^2 \ll \varepsilon_1$ in \eqref{ch1_spat_oprtr_2}, the solution has regular boundary layers near all the four edges and also the corner layers near every corner of $\varOmega$. The regular boundary layers have width of $\displaystyle O(\sqrt{\eps_1})$ \cite{o2006numerical}. \textcolor{black}{To demonstrate the uniform convergence of numerical methods for SPPs on layer-adapted meshes, it is essential to decompose the solution into a sum of regular and layer components}.

The following result holds for problem \eqref{elliptic_problem} when the operator is defined as in \eqref{ch1_spat_oprtr_2} \cite{o2006numerical}.
\begin{theorem}
    \textcolor{black}{Assume that the solution $\mathfrak{v}$ of \eqref{elliptic_problem} has the following decomposition}
    \begin{equation}
	\mathfrak{v} = \mathscr{S} + \mathscr{E}_L + \mathscr{E}_R + \mathscr{E}_T + \mathscr{E}_B + \mathscr{E}_{LT} + \mathscr{E}_{RB} + \mathscr{E}_{RT} + \mathscr{E}_{LB},
\end{equation}
where $\mathscr{S}$ is the smooth part while $\mathscr{E}_L ,\, \mathscr{E}_R ,\, \mathscr{E}_B ,\, \mathscr{E}_T$ are the boundary layer parts \textcolor{black}{near the domain boundary} $\mathcal{D}_{10},\, \mathcal{D}_{11},\, \mathcal{D}_{20},\, \mathcal{D}_{21}$ respectively and $\mathscr{E}_{LB} ,\,  \mathscr{E}_{LT} ,\,  \mathscr{E}_{RB} ,\,  \mathscr{E}_{RT}$ are the corresponding corner layer parts appearing at the corners $(0,0),\, (1,0),\, (0,1),\, (1,1)$, where
\begin{align*}
\mathcal{D}_{10} = \big\{(0,y):\,\, 0 \leq  y \leq 1\big\},\quad \mathcal{D}_{11} = \big\{(1,y):\,\, 0 \leq  y \leq 1\big\},\\[6 pt]
\mathcal{D}_{20} = \big\{(x,0):\,\, 0 \leq  x \leq 1\big\},\quad \mathcal{D}_{21} = \big\{(x,1):\,\, 0 \leq  x \leq 1\big\}.
\end{align*}
\textcolor{black}{
Let $\beta = \min \big\{ \beta_1, \beta_2 \big\} $ and $\displaystyle \chi < \min_{\overline\Omega} \bigg\{\dfrac{\mathfrak{r}(x,y)}{2 \mathfrak{d}_1(x,y)},\,\ \dfrac{\mathfrak{r}(x,y)}{2 \mathfrak{d}_2(x,y)} \bigg\}$.
The various components of the solution $\mathfrak{v}$ to the problem \eqref{elliptic_problem} satisfy the explicit bounds given below. If the compatibility conditions in \eqref{ch1_compatible} are met, then the boundary layer components adhere to the following bounds:
\begin{eqnarray*}
\left\Vert\dfrac{\partial^{l_{1}+l_{2}}\mathscr{S}}{\partial x^{l_1} \partial y^{l_2}}\right\Vert &\leq& C(1 + \varepsilon^{1 - (l_1 + l_2)/2}), \quad 0 \leq l_1 + l_2 \leq 3, \\ [2pt]
|\mathscr{E}_{L}(x,y)| &\leq& C\exp\Big(-\sqrt{\chi \beta/\varepsilon_1} \, x\Big), \\ [2pt]
|\mathscr{E}_{B}(x,y)| &\leq&  C\exp\Big(-\sqrt{\chi \beta/\varepsilon_1} \, y\Big), \\ [2pt]
|\mathscr{E}_{R}(x,y)| &\leq&  C\exp\Big(-\sqrt{\chi \beta/\varepsilon_1} \, (1-x)\Big), \\ [2pt]
|\mathscr{E}_{T}(x,y)| &\leq&  C\exp\Big(-\sqrt{\chi \beta/\varepsilon_1} \, (1-y)\Big), \\ [2pt]
|\mathscr{E}_{LB}(x,y)| &\leq&  C\exp\Big(-\sqrt{\chi \beta/\varepsilon_1} \, x\Big)\exp\Big(-\sqrt{\chi \beta/\varepsilon_1} \, y\Big), \\ [2pt]
|\mathscr{E}_{LT}(x,y)| &\leq& C\exp\Big(-\sqrt{\chi \beta/\varepsilon_1} \, x\Big)\exp\Big(-\sqrt{\chi \beta/\varepsilon_1} \, (1-y)\Big), \\ [2pt]
|\mathscr{E}_{RB}(x,y)| &\leq&  C\exp\Big(-\sqrt{\chi \beta/\varepsilon_1} \, (1-x)\Big)\exp\Big(-\sqrt{\chi \beta/\varepsilon_1} \, y\Big),\\ [2pt]
|\mathscr{E}_{RT}(x,y)| &\leq&  C\exp\Big(-\sqrt{\chi \beta/\varepsilon_1} \, (1-x)\Big)\exp\Big(-\sqrt{\chi \beta/\varepsilon_1} \, (1-y)\Big), \\ [2pt]
\left\Vert \dfrac{\partial^{l_1}\mathscr{E}_L}{\partial y^{l_1}}\right\Vert &\leq& C(1 + \varepsilon_1^{(1-l_1)/2}), \qquad \, \left\Vert \dfrac{\partial^{l_1}\mathscr{E}_R}{\partial y^{l_1}}\right\Vert \leq C(1 + \varepsilon_1^{(1-l_1)/2}),\, 1 \leq l_1 \leq 3, \\ [2pt]
\left\Vert \dfrac{\partial^{l_1}\mathscr{E}_B}{\partial x^{l_1}}\right\Vert &\leq& C(1 + \varepsilon_1^{(1-l_1)/2}), \qquad \, \left\Vert \dfrac{\partial^{l_1}\mathscr{E}_T}{\partial x^{l_1}}\right\Vert \leq C(1 + \varepsilon_1^{(1-l_1)/2}),\, 1 \leq l_1 \leq 3.
\end{eqnarray*}
We also have the following bound for each layer component,
\begin{equation*}
    \left\Vert\dfrac{\partial^{l_{1}+l_{2}}\mathscr{E}}{\partial x^{l_1} \partial y^{l_2}}\right\Vert \leq C\varepsilon^{-(l_1 + l_2)/2}, \quad 1 \leq l_1 + l_2 \leq 3.
\end{equation*}
}
\end{theorem}
Similar result holds for the case when the operator is defined as in \eqref{ch1_spat_oprtr} and can be found in \cite{teofanov2007elliptic}.

\subsubsection{ 2D Parabolic initial-boundary-value problem }

Consider the following class of 2D singularly perturbed parabolic PDEs with two parameters:
\begin{equation}\label{ch3_model_prob}
\left\{
\begin{array}{lll}
\mathfrak{v}_t  + \mbox{{L}}_{\varepsilon_1,\varepsilon_2}  \mathfrak{v} = h(x,y,t), &(x,y)\in \varOmega = (0,1) \times (0,1),\,\, t \in \varOmega_t,\\ [6pt]
\mathfrak{v}(x,y,t) = 0, &(x,y) \in \partial \varOmega = \overline{\varOmega} \setminus \varOmega, \,\, t \in \overline{\varOmega}_t,\\ [6pt]
\mathfrak{v}(x,y,0) = \mathfrak{v}_{0}(x,y), & (x,y) \in \varOmega,
\end{array} \right.
\end{equation}
where $\mbox{{L}}_{\varepsilon_1,\varepsilon_2}$ is same as defined for the elliptic case. Here $\varOmega_t = (0, \mathscr{T}]$ and $h$ satisfies the following compatibility restriction:
\begin{equation}
    h(x,0,t) = h(x,1,t) = h(0,y,t) = h(1,y,t) = 0.
\end{equation}

\section{Modified Physics Informed Neural Network }\label{modi_pinn}

\subsection{Fully connected neural network: Mathematical background}\label{fully_connected}

Let \(\mathscr{F}: \mathbb{R}^{m} \rightarrow \mathbb{R}^{m_{out}}\) denotes a feed-forward neural network with \(s-1\) hidden layers and \(n_l\) neurons in the \(l\)-th layer. In the \(l\)-th layer, biases are represented by the vector \(\Tilde{\textbf{b}}^l \in \mathbb{R}^{n_l}\), and weights by the matrix \({W}^{l} \in \mathbb{R}^{n_l \times n_{l-1}}\). Then, we can define $\mathscr{F}$ as follows:
For $x \in \mathbb{R}^{n_0}$ and $\mathscr{F}^{l}(x) \in \mathbb{R}^{n_l}$:
\begin{equation}
    \mathscr{F}^{l}(x) = {W}^{l}\beta(\mathscr{F}^{l-1}(x)) + \Tilde{\textbf{b}}^l \in \mathbb{R}^{l}, \,\,\, \mbox{for}\,\, 2 \leq l\leq s
\end{equation}
and
\[
 \mathscr{F}^{1}(x) = {W}^{1}(x) + \Tilde{\textbf{b}}^1\, ,
\]
where $\beta$ is a nonlinear activation function. Now the overall output can be written as
\[
\mathfrak{v}_{\theta} = \mathscr{F}^{s} \circ \textcolor{black}{\beta} \circ \mathscr{F}^{s-1}\circ \ldots \circ \textcolor{black}{\beta} \circ \mathscr{F}^{1}(x),
\]
where $\theta = \{{W}^{l},\Tilde{\textbf{b}}^l\},\,\, l = 1,2,\ldots,s$. A basic feed-forward neural network is shown in Figure \ref{FFRN}.

\textbf{\begin{figure}[ht]
  \centering
  \includegraphics[width=0.8\textwidth]{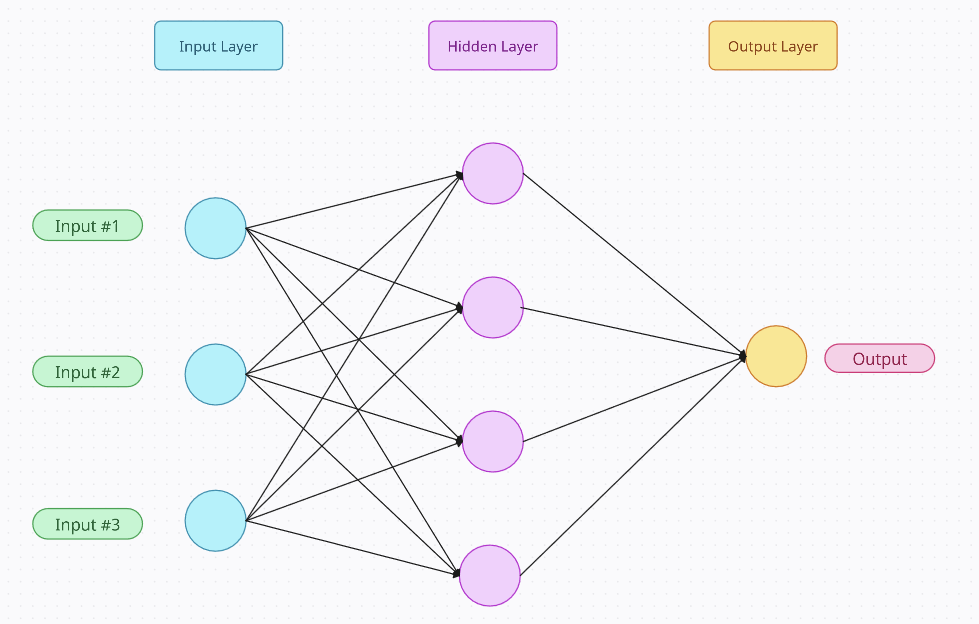}
  \caption{A basic feed-forward neural network.}
  \label{FFRN}
\end{figure}}

\subsection{Parameter Asymptotic PINNs (PA-PINNs) Methodology \textcolor{black}{\cite{cao2023physics}}}
PINNs are a form of unsupervised machine learning capable of estimating solutions to PDEs without requiring a labeled training dataset. They convert the task of solving the governing equations into an optimization problem aimed at finding a set of parameters that minimizes a loss function. This function is built from the residuals arising from the governing equations, boundary conditions, and initial conditions.

Consider a differential equation in its general form as follows:
\begin{equation}\label{model_pb}
\left\{
\begin{array}{ll}
\mathfrak{D}[\mathfrak{v}(x);\zeta] = \textcolor{black}{h(x)} , \quad  x \in \varOmega, \\[6pt]
\mathfrak{BC}[\mathfrak{v}(x)] = \textcolor{black}{B(x)}, \quad x \in \partial \varOmega,
\end{array}\right.
\end{equation}
where $\varOmega$ is a bounded domain in $\mathbb{R}^{n}$, $n$ denotes the dimension, with  $\mathfrak{D}$ and $\mathfrak{BC}$ being the differential operator and boundary operator respectively. Here $\mathfrak{v}$ denotes the exact solution, \textcolor{black}{$h$} represents the source function, $\zeta$ are the parameters and \textcolor{black}{$B$} is the specified boundary function.
Now, the solution $\mathfrak{v}(x)$ is approximated using PINNs by training a neural network $N(x;\theta)$, where $\theta$ represents the network's hyperparameters. The objective is to minimize the following loss function, which will guide the training process:
\begin{equation}\label{loss_pinn}
    \mathscr{L}^{total}(\theta) = \mathscr{L}^{o}(\theta) + \mathscr{L}^{b}(\theta),
\end{equation}
where  
\begin{eqnarray}
\mathscr{L}^{o}(\theta) &=& \dfrac{1}{\eta_{o}} \sum_{i=1}^{\eta_{o}}\|\mathfrak{D}[N(x_{i}^{o};\theta)]-\textcolor{black}{h}(x_{i}^{o})\|^{2}, \\[6pt]
\mathscr{L}^{b}(\theta) &=& \dfrac{1}{\eta_{b}} \sum_{i=1}^{\eta_{b}}\|\mathfrak{BC}[N(x_{i}^{b};\theta)] - \textcolor{black}{B}(x_{i}^{b})\|^{2}.
\end{eqnarray}
The terms $\mathscr{L}^{o}(\theta)$ and $\mathscr{L}^{b}(\theta)$ represent the ``operator loss" and ``boundary loss," respectively. Here, $\eta_o$ and $\eta_b$ denote the number of collocation points sampled in the interior $\varOmega$ and at the boundary $\partial \varOmega$, respectively. Also $ \{x_{i}^{o}\}$ and $\{x_{i}^{b}\}$ denote the sets of training points taken from the interior and boundary of the domain $\varOmega$, respectively.

When training a PINNs network, several key considerations must be addressed. It is essential to ensure an adequate number of data points are sampled to facilitate the network in learning a cohesive solution across the entire domain. To prevent the network from overfitting one should consider applying regularization techniques \textcolor{black}{\cite{bajaj2023recipes}}. Furthermore, the loss function needs to be differentiated with respect to the hyperparameters $\theta$ of the neural network. Optimization techniques such as gradient descent are then applied to optimize \( \theta \). It has been documented in the literature that standard PINNs often exhibit lower convergence and accuracy for solving SPPs and may fail to produce satisfactory results \textcolor{black}{\cite{karniadakis2021physics}}.

The primary concept behind the PA-PINNs methodology involves initially approximating the smooth part of the solution using large parameter values. Subsequently, the boundary layer solution is approximated by gradually decreasing the perturbation parameter. This approach eliminates the need for apriori information about the position of the boundary layer. Following each parameter reduction, a PINN is trained to converge using optimal network parameters. An adaptive algorithm is then employed to refine the distribution of training points based on the loss function, updating the training set accordingly. This adaptive algorithm allows for selecting more points from the boundary layer region, thereby enhancing both the accuracy and convergence of the networks.

For illustration, we can rewrite problem \eqref{elliptic_problem} as follows
\begin{equation}
    \left\{\begin{array}{ll}
\mathfrak{L}_{\varepsilon_{1}^{j},\varepsilon_{2}^{j}}\mathfrak{v} := \mbox{L}_{\varepsilon_{1}^{j},\varepsilon_{2}^{j}}\mathfrak{v}({\bf x}) - h({\bf x}) = 0, \quad \mbox{for} \, {\bf x} \in \varOmega,\\ [6pt]
\mathfrak{B}\mathfrak{v} := \mathcal{B}\mathfrak{v}({\bf x}) = 0,\quad \mbox{for} \, {\bf x} \in \partial\varOmega. \quad \\ [6pt]
\end{array}\right.
\end{equation}
Then the associated residual loss function will be
\begin{equation}
    \mathcal{L}(\eta) = \vartheta_{1} \times \dfrac{1}{\textcolor{black}{\eta\strut_{o}}} \sum_{j=1}^{\textcolor{black}{\eta\strut_{o}}}\|\mbox{L}_{\varepsilon_{1},\varepsilon_{2}}\mathfrak{v}(\textcolor{black}{{\bf x}_{j}^{o}}) - h(\textcolor{black}{{\bf x}_{j}^{o}})\|^{2} + \vartheta_{2}\times\dfrac{1}{\eta\strut_{b}} \sum_{j=1}^{\eta\strut_{b}}\|\mathcal{B}_{}\mathfrak{v}({\bf x}_{j}^{b})\|^{2},
\end{equation}
where $\vartheta_1$ and $\vartheta_2$ are weights, \textcolor{black}{$\{{\bf x}_{j}^{o}\}$} are the collocation points from inside the domain $\varOmega$ and $\textcolor{black}{\{{\bf x}_{j}^{b}\}} \in  \partial \varOmega$  are the collocation points from the initial/boundary location. Here $\eta\strut_{t} = \textcolor{black}{\eta\strut_{o}}+\eta\strut_{b}$ is the training set cardinality where $\textcolor{black}{\{{\bf x}\}} = \textcolor{black}{\{{\bf x}^{o}\}} \cup \textcolor{black}{\{{\bf x}^{b}\}}$. Also, $\varepsilon_{1}^{j}$ and $\varepsilon_{2}^{j}$ denote the values of the perturbation parameters.
Choosing an update rule for $\varepsilon_{i}^{j}$ as follows
\begin{equation}\label{linear_update}
\varepsilon_{i}^{j+1} = \varepsilon_{i}^{j} - \delta_{\varepsilon_{i}}, \quad  i = 1,2 ;\,\, j = 0,1,\ldots,N,
\end{equation}
where $\delta_{\varepsilon_{i}} = (\varepsilon_{i}^{0} - \varepsilon_{i}^{N})/(N+1)$.

\textcolor{black}{
\begin{remark}
    As can be seen that we have considered linear update rule for the perturbation parameters. One can also choose some nonlinear update rule as well. Suppose in our case $\varepsilon_{1}^{N}$ and $\varepsilon_{2}^{M}$ be the final parameter values for $\varepsilon_1$ and $\varepsilon_2$ respectively. Then one can choose the update as follows:
    \begin{equation*}
        \varepsilon_1^{k} = \varepsilon_{1}^{0}\,c_{1}^{k}, \quad k = 1,2,\ldots,N, \quad N =\Big[ \log_{c_1}\left(\dfrac{\varepsilon_{1}^{N}}{\varepsilon_1^{0}}\right) \Big]
    \end{equation*}
    and
    \begin{equation*}
        \varepsilon_2^{p} = \varepsilon_{2}^{0}\,c_{2}^{p}, \quad p = 1,2,\ldots,M, \quad M =\Big[ \log_{c_2}\left(\dfrac{\varepsilon_{2}^{M}}{\varepsilon_2^{0}}\right) \Big],
    \end{equation*}
    where $\varepsilon_{1}^{0}$ and $\varepsilon_{2}^{0}$ are the initial values of the perturbation parameters. The parameter $c_1 = 0.71$  and $c_2 = 0.72$ are used to obtain parameter partition. Other values of $c_1$ and $c_2$ can also be used.
\end{remark}
}

\subsection{Algorithm for PA-PINNs}
The PA-PINN algorithm for two-parameter problem is illustrated in \textbf{Algorithm} \ref{alg:1}. 
Figure \ref{work_flow1} provides a visual representation of the first half of the algorithm, illustrating how $\varepsilon_{1}$ progresses to $\varepsilon_{1}^{end}$. A similar process is followed for the parameter $\varepsilon_{2}$.

\begin{algorithm}[H]
\label{alg:1}
\noindent\rule[0.5ex]{\linewidth}{1pt}
{\bf Algorithm 1: PA-PINN Algorithm }\\
\noindent\rule[0.5ex]{\linewidth}{1pt}
Initialize $\textcolor{black}{\vartheta^{0}}, \mathfrak{d}^{0}$, initial training data $\{{\bf x}\}_{0}$, $\varepsilon_{1}^{0}$ and $\varepsilon_{2}^{0}$ and $j=0$. \newline
\While{$\varepsilon^{j}_{1} \geq \varepsilon^{end}_{1}$ and keeping $\varepsilon_{2}$ fixed at $\varepsilon_{2}^{0}$} {
\For{$i = 1,2,\ldots,$} {
Solve for $\Bar{u}^{j,i}$ \newline
$L_{\varepsilon_{1}^{j},\varepsilon_{2}^{0}}u := \mathcal{L}_{\varepsilon_{1}^{j},\varepsilon_{2}^{0}}u({\bf x}) - h({\bf x}) = 0, \,\, \mbox{for} \, {\bf x} \in \Omega,$ \newline
 $\mathfrak{B}u := \mathcal{B}u({\bf x}) = 0, \,\, \mbox{for} \, {\bf x} \in \partial\Omega,$ \newline
if $\|\Bar{u}^{j,i}-\Bar{u}^{j,i-1}\|/\|\Bar{u}^{j,i}\| < tol$ then $\Bar{u}^{i} = \Bar{u}^{j,N}$ \newline
}
Adaptive refinement based on residual to obtain new testing/collocation points $\delta \{{\bf x}\}$\newline
Updating parameter $\varepsilon_{1}^{j+1} = \varepsilon_{1}^{j} - \delta_{\varepsilon_1}$ or $\varepsilon_1^{j+1} = \varepsilon_{1}^{j}\,c_{1}^{j+1}$, collocation points  $\{{\bf x}\}_{j+1} = \{{\bf x}\}_{j} + \delta \{{\bf x}\}$\newline
\textcolor{black}{Use updated weights as initial weights for the next loop.}\newline
}
Now \newline
\While{$\varepsilon^{j}_{2} \geq \varepsilon^{end}_{2}$ and $\varepsilon_{1}$ reached at $\varepsilon_{1}^{end}$ from the earlier loop} {
\For{$i = 1,2,\ldots,$} {
Solve for $\Bar{u}^{j,i}$ \newline
$L_{\varepsilon_{1}^{end},\varepsilon_{2}^{j}}u := \mathcal{L}_{\varepsilon_{1}^{end},\varepsilon_{2}^{j}}u({\bf x}) - h({\bf x}) = 0, \,\, \mbox{for} \, {\bf x} \in \Omega,$ \newline
 $\mathfrak{B}u := \mathcal{B}u({\bf x}) = 0, \,\, \mbox{for} \, {\bf x} \in \partial\Omega,$ \newline
if $\|\Bar{u}^{j,i}-\Bar{u}^{j,i-1}\|/\|\Bar{u}^{j,i}\| < tol$ then $\Bar{u}^{i} = \Bar{u}^{j,N}$ \newline
}
Adaptive refinement based on residual to obtain new testing/collocation points $\delta \{{\bf x}\}$\newline
Updating parameter $\varepsilon_{2}^{j+1} = \varepsilon_{2}^{j} - \delta_{\varepsilon_2}$ or $\varepsilon_2^{j+1} = \varepsilon_{2}^{j}\,c_{2}^{j+1}$, collocation points  $\{{\bf x}\}_{j+1} = \{{\bf x}\}_{j} + \delta \{{\bf x}\}$\newline
}
\noindent\rule[0.5ex]{\linewidth}{1pt}
\end{algorithm}

\textbf{\begin{figure}[ht]
  \centering
  \includegraphics[width=1.0\textwidth]{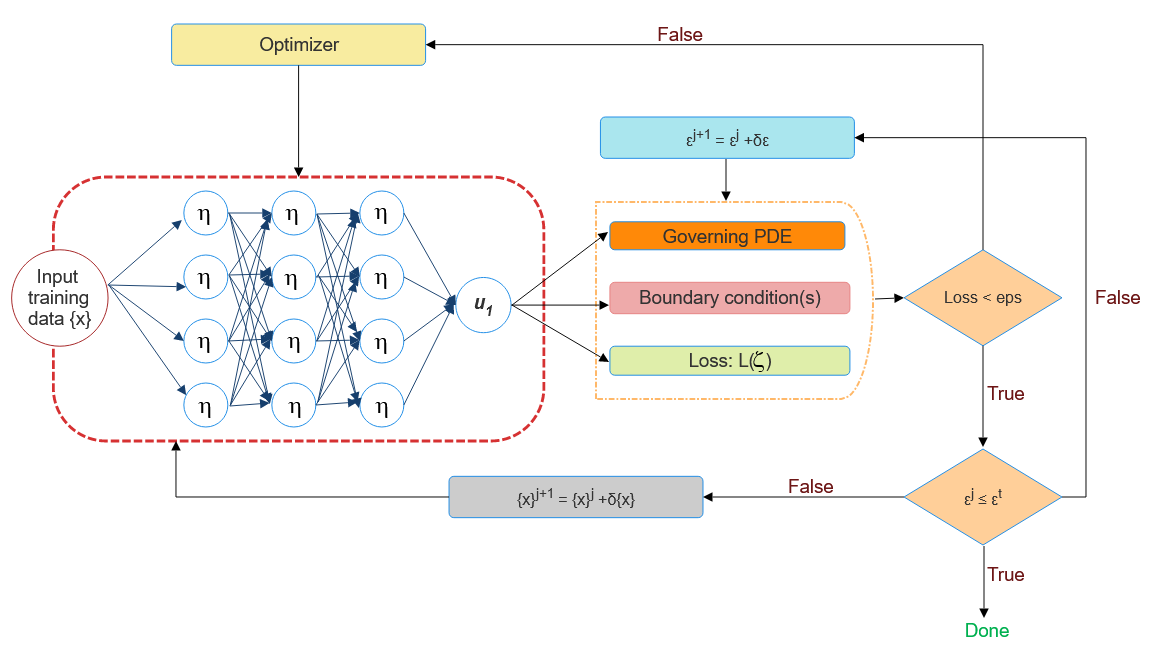}
  \caption{PA-PINNs Algorithm flowchart.}
  \label{work_flow1}
\end{figure}}

\section{Numerical Experiments} \label{num_ex}
\textcolor{black}{In this section, we illustrate the effectiveness of PA-PINNs in addressing two-parameter SPPs through a series of examples.}	 
For implementation we have used fully connected feed-forward neural network. For 1D problems we have taken $\eta_{o}$ as $10^{3}$ and $\eta_{b}$ as $\textcolor{black}{256}$, where $\eta_{o}$ and $\eta_{b}$ denote the number of collocation points from the interior and boundary of the domain respectively taken for the training. For time-independent problems our network has 8 hidden layers with 20 neurons in each layer and for time-dependent problems network consists of 2 hidden layers with 20 neurons each. L-BFGS optimizer has been used and activation function is taken as $\tanh$. The Latin hypercube sampling method is applied to gather points inside the interior domain and along the boundaries. For 2D problems we have taken $\eta_{o}$ as $10^{4}$ and $\eta_{b}$ as $10^3$. We ran our algorithm for 4000 iterations. For all the examples we have choosen $(\varepsilon_{1}^{0}, \varepsilon_{2}^{0}) = (0.3, 0.1)$. All other network parameters are identical as was for 1D problems. We will use a linear update rule unless stated otherwise.

We are using the relative $L^{2}$ form of error and is given as
\[
\displaystyle\mathscr{E}_{2} (\mathfrak{v}, \hat{\mathfrak{v}}) = \sqrt{\dfrac{\sum_{j=0}^{N - 1} (\mathfrak{v}_j - \hat{\mathfrak{v}}_j)^2}{\sum_{j=0}^{N - 1} \mathfrak{v}_j^2}},
\]
where $\mathfrak{v}$ and $\hat{\mathfrak{v}}$ are the exact and neural network solutions respectively. For Example \ref{ex2D_time} we have also provided the $L^{\infty}$ error, defined as
\[
\mathscr{E}_{\infty}(\mathfrak{v}, \hat{\mathfrak{v}}) = \max_{0 \leq j \leq N-1}\{|\mathfrak{v}_{j} - \hat{\mathfrak{v}}_{j}|\}.
\]
\textcolor{black}{In Tables \ref{table:1}--\ref{table:7}, we present the errors corresponding to linear update rule \eqref{linear_update} and a comparison with standard PINNs}.

\begin{figure}[ht]
  \centering
  \subfloat[$(\varepsilon_1,\varepsilon_2) = (10^{-2},10^{-3})$]{\includegraphics[width=.33\linewidth]{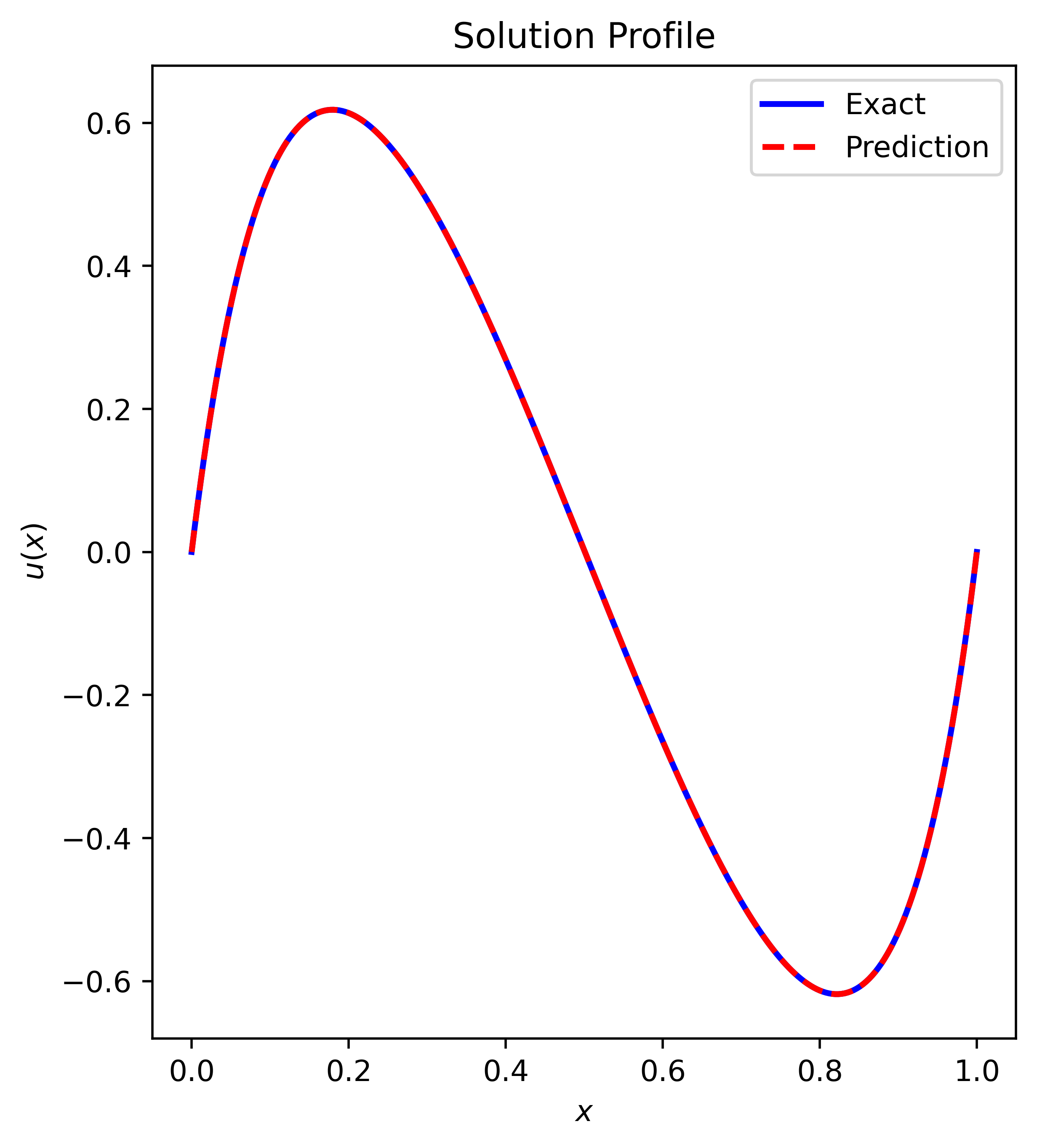}}\hfill
  \subfloat[$(\varepsilon_1,\varepsilon_2) = (10^{-3},10^{-4})$] {\includegraphics[width=.33\linewidth]{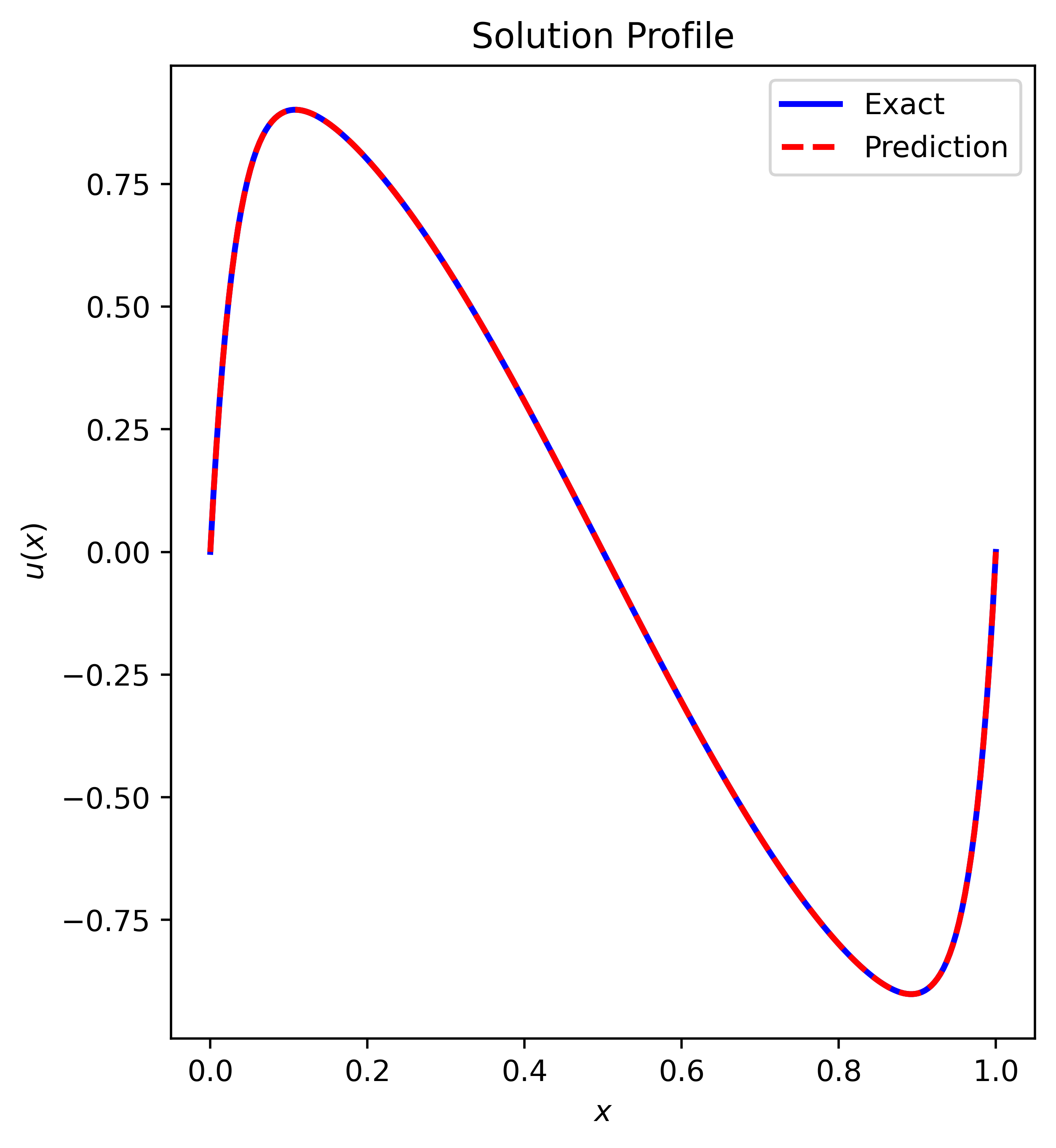}}
  \subfloat[$(\varepsilon_1,\varepsilon_2) = (10^{-4},10^{-5})$] {\includegraphics[width=.33\linewidth]{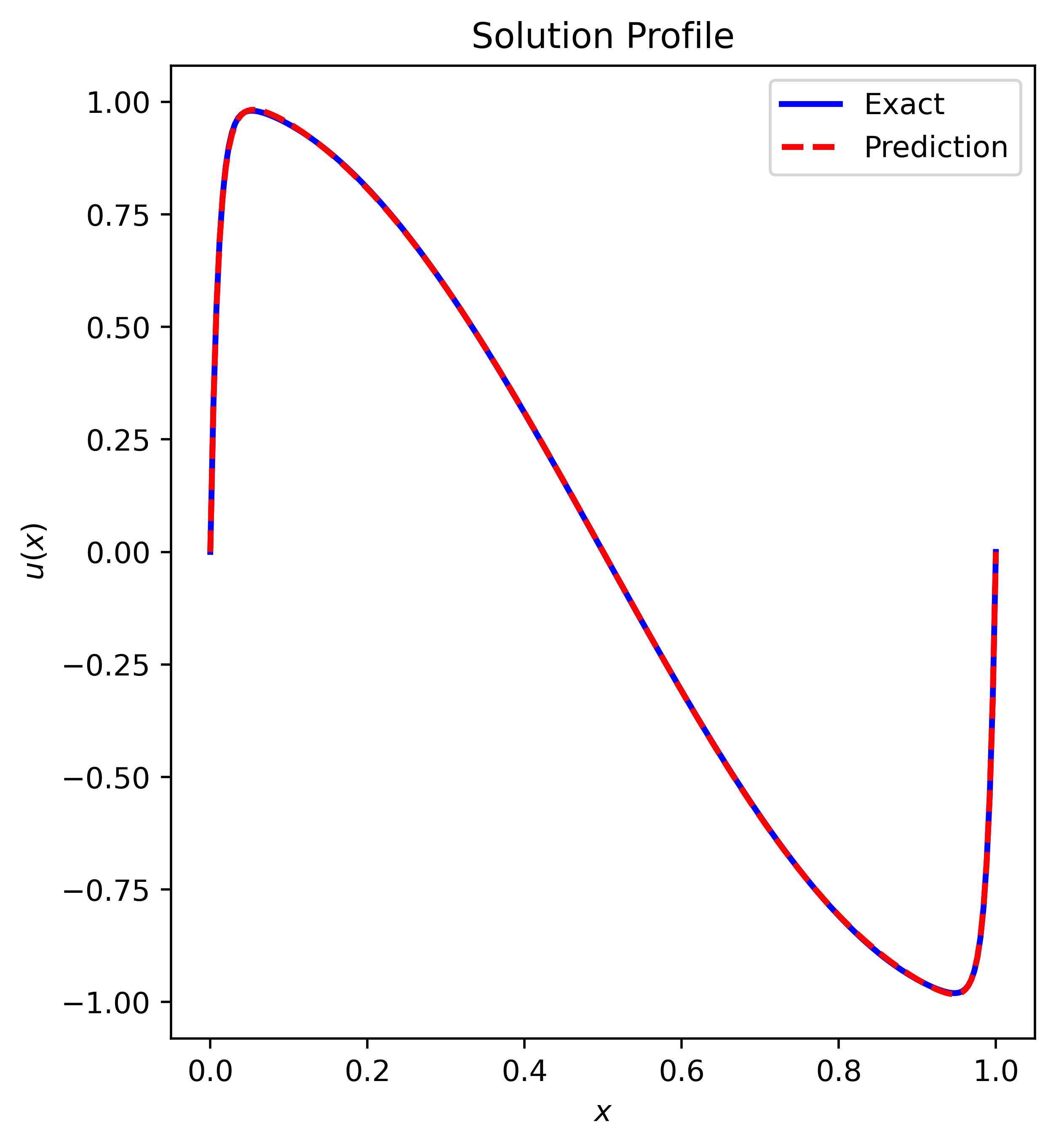}}
  \caption{Comparison between Exact and PA-PINNs solution for different values of perturbation parameter $\varepsilon_1, \varepsilon_2$ for Example \ref{ch2exmp1}.}\label{fig1_ex1}
\end{figure}

\begin{example}\label{ch2exmp1}
		Consider the following test problem 
		\begin{equation}
		\left\{
		\begin{array}{ll}
		-\varepsilon_1 \mathfrak{v}''+\varepsilon_2 \mathfrak{v}'+ \mathfrak{v}=\cos\pi x,
		\quad x\, \in (0,1),  &\\[6pt]
		\qquad \mathfrak{v}(0) = 0 = \mathfrak{v}(1),
		\end{array} \right.
		\end{equation}
    whose exact solution is represented by
	\[
	\mathfrak{v}(x)=\mathfrak{A}_1\cos \pi x + \mathfrak{A}_2\exp(-\mu_Lx) + \mathcal{B}_1 \sin \pi x + \mathcal{B}_2\exp(-\mu_R(1-x)),
	\]
	where
	\[
	 \mathfrak{A}_1=\frac{\varepsilon_1\pi^2+1}{\varepsilon_2^2\pi^2+(\varepsilon_1\pi^2+1)^2}, \quad \mathcal{B}_1=\frac{\varepsilon_2\pi}{\varepsilon_2^2\pi^2+(\varepsilon_1\pi^2+1)^2},
	\]
	\[
	 \mathfrak{A}_2=-\mathfrak{A}_1\frac{1+\exp(-\mu_R)}{1-\exp(-(\mu_L+\mu_R))}, \quad \mathcal{B}_2=\mathfrak{A}_1\frac{1+\exp(-\mu_L)}{1-\exp(-(\mu_L+\mu_R))}, \quad \mu_{L,R} = \frac{\mp \varepsilon_2 + \sqrt{\varepsilon_2^2+4\varepsilon_1}}{2\varepsilon_1}.
	\]
\end{example}
We ran the code for $1000$ iterations and the errors obtained are presented in Table \ref{table:1}. The table shows that PA-PINNs perform significantly better than traditional PINNs. Solution plots for various parameter values are provided in the Figure \ref{fig1_ex1}. From these plots, it is evident that as the parameter values decrease, the sharpness of the boundary layers increases, and PA-PINN effectively captures the solution profile with accuracy.

\begin{table}[!h]
\centering
\caption{Error comparison for various parameter values for Example \ref{ch2exmp1}.}
\label{table:1}
\begin{tabular}{|c|cc|c|c|}
\hline \text { Method } & $\varepsilon_{1}$ & $\varepsilon_2$ & $\mathscr{E}_{2}$\\
\hline & $10^{-2}$ & $10^{-3}$ & \textcolor{black}{1.661e-03}\\
\textcolor{black}{\text { PINN }} & $10^{-3}$ & $10^{-4}$ & \textcolor{black}{1.287e-02} \\
& $10^{-4}$ & $10^{-5}$ & \textcolor{black}{7.905e-02}\\
& $10^{-5}$ & $10^{-6}$ & \textcolor{black}{1.357e-01}\\
\hline & $10^{-2}$ & $10^{-3}$ & 2.446e-05\\
\text { PA-PINN }& $10^{-3}$ & $10^{-4}$ & 1.548e-04\\
 & $10^{-4}$ & $10^{-5}$ & 1.782e-03\\
& $10^{-5}$ & $10^{-6}$ & 1.273e-02\\
\hline
\end{tabular}
\end{table}

\begin{example}\label{ch2exmp2}
    Consider the following two-parameter problem
    \begin{equation}
		\left\{
		\begin{array}{ll}
		-\varepsilon_1 \mathfrak{v}''-\varepsilon_2 \mathfrak{v}'+ \mathfrak{v}= \exp(1-x)
		\quad x\, \in (0,1),  &\\[6pt]
		\qquad \mathfrak{v}(0) = 0 = \mathfrak{v}(1),
		\end{array} \right.
		\end{equation}
  whose exact solution is
  \begin{eqnarray*}
          \mathfrak{v}(x) &=& \dfrac{1}{\varepsilon_1 - \varepsilon_2 - 1}\Bigg[\dfrac{\exp(1) - \exp(\mathfrak{m}_{1})}{1-\exp(-\sqrt{\varepsilon_2^{2} + 4\varepsilon_1}/\varepsilon_1)} \times \exp(-\mathfrak{m}_{2}x) - \exp(1-x) \\ [4pt] 
          && \qquad \qquad \qquad + \dfrac{1-\exp{(1-\mathfrak{m}_{2})}}{1-\exp(-\sqrt{\varepsilon_2^{2} + 4\varepsilon_1}/\varepsilon_1)}\times \exp(\mathfrak{m}_{1}(1-x))\Bigg],
  \end{eqnarray*}
  where $\mathfrak{m}_{1,2} = \dfrac{\varepsilon_2 \mp \sqrt{\varepsilon_2^{2} + 4\varepsilon_1}}{2\varepsilon_1}$.
\end{example}

\begin{figure}[ht]
  \centering
  \subfloat[$(\varepsilon_1,\varepsilon_2) = (10^{-2},10^{-3})$]{\includegraphics[width=.33\linewidth]{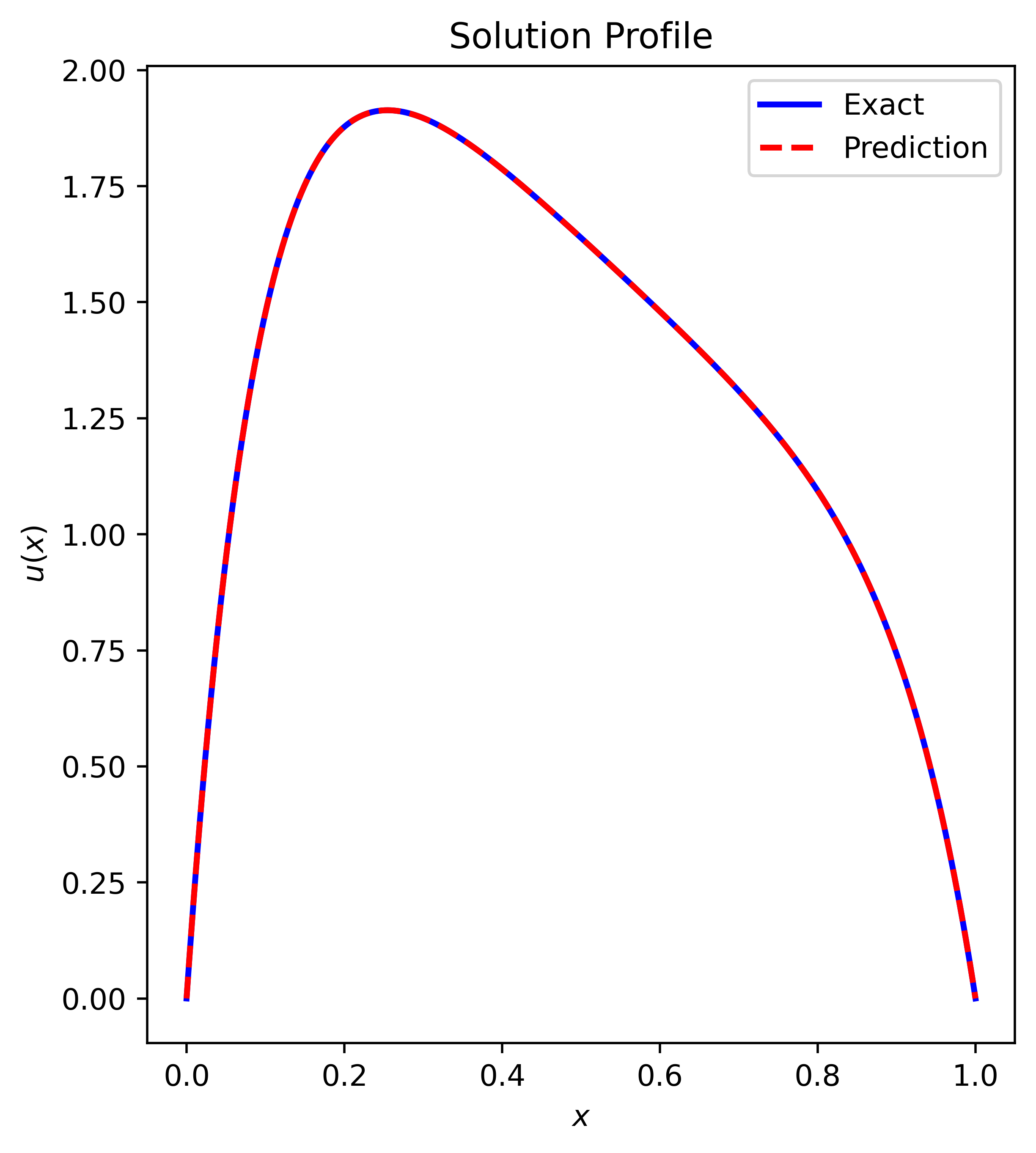}}\hfill
  \subfloat[$(\varepsilon_1,\varepsilon_2) = (10^{-3},10^{-4})$] {\includegraphics[width=.33\linewidth]{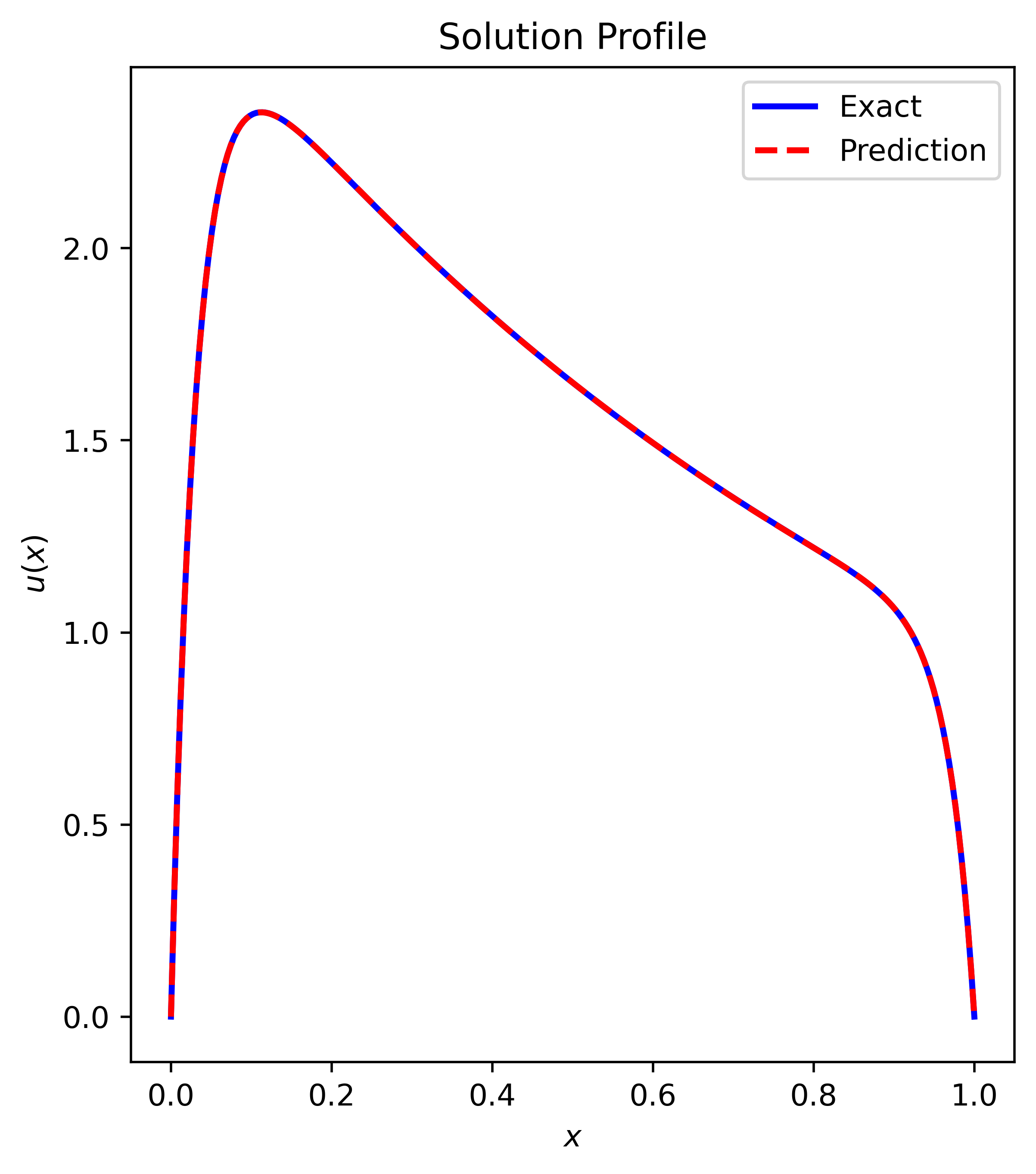}}
  \subfloat[$(\varepsilon_1,\varepsilon_2) = (10^{-4},10^{-5})$] {\includegraphics[width=.33\linewidth]{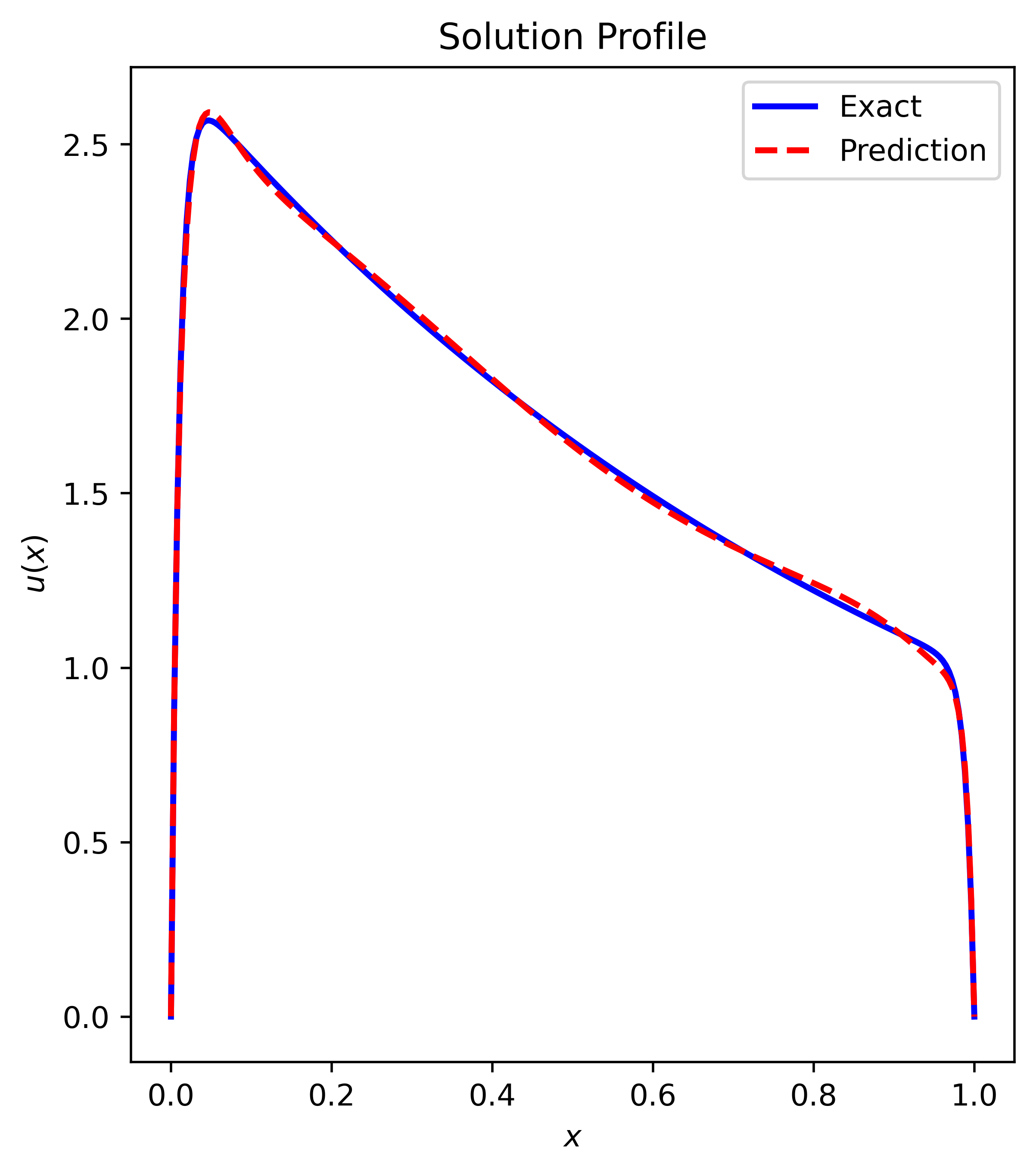}}
  \caption{Comparison between Exact and PINNs solution for different values of perturbation parameter $\varepsilon_1, \varepsilon_2$ for Example \ref{ch2exmp2}.}\label{fig1_ex4}
\end{figure}

We executed the code for 1000 iterations and documented the errors in Table \ref{table:3}, which compares the result of PA-PINNs with those of traditional PINNs, showing that PA-PINNs significantly outperforms them. Figure \ref{fig1_ex4} presents solution plots for different parameter values, showing that as the parameter values decrease, the boundary layers become sharper, and PA-PINN accurately captures the solution profile.

\begin{table}[!h]
\centering
\caption{Error comparison for various parameter values for Example \ref{ch2exmp2}.}
\label{table:3}
\begin{tabular}{|c|cc|c|c|}
\hline \text { Method } & $\varepsilon_{1}$ & $\varepsilon_2$ & $\mathscr{E}_{2}$\\
\hline & $10^{-1}$ & $10^{-2}$ & \textcolor{black}{1.911e-04}\\
\textcolor{black}{\text { PINN }} & $10^{-2}$ & $10^{-3}$ & \textcolor{black}{2.312e-03} \\
& $10^{-3}$ & $10^{-4}$ & \textcolor{black}{1.790e-01}\\
& $10^{-4}$ & $10^{-5}$ & \textcolor{black}{4.577e-02}\\
\hline & $10^{-1}$ & $10^{-2}$ & 4.378e-04 \\
\text { PA-PINN }& $10^{-2}$ & $10^{-3}$ & 2.859e-06 \\
 & $10^{-3}$ & $10^{-4}$ & 2.295e-05 \\
& $10^{-4}$ & $10^{-5}$ & 8.580e-03\\
\hline
\end{tabular}
\end{table}

\begin{example}\label{ch2exmp3}
		Consider the following parabolic initial-boundary-value problem, which is the time-dependent version of  Example \ref{ch2exmp1}:
		\begin{equation}\label{}
		\left\{
		\begin{array}{ll}
		\dfrac{\partial \mathfrak{v}}{\partial t}-\varepsilon_1 \dfrac{\partial^2\mathfrak{v}}{\partial x^2} + \varepsilon_2 \dfrac{\partial\mathfrak{v}}{\partial x} + (1-\exp(-t))\mathfrak{v}= (1-\exp(-t))\cos\pi x,
		\quad x\, \in (0,1), \, 0\leq t\leq \mathscr{T},  &\\[10pt]
		\mathfrak{v}(x,0) = 0, \, \forall \, x \in (0,1), \quad \mathfrak{v}(0,t) = 0 = \mathfrak{v}(1,t),
		\end{array} \right.
		\end{equation}
whose exact solution is represented by
\[
\mathfrak{v}(x)=\Big(1-\exp(-t)\Big)\times\left(\mathfrak{A}_1\cos \pi x + \mathfrak{A}_2\exp(-\mu_Lx)+  \mathcal{B}_1 \sin \pi x + \mathcal{B}_2\exp(-\mu_R(1-x))\right),
\]
where $\mathfrak{A}_1, \mathfrak{A}_2, \mathfrak{B}_1, \mathfrak{B}_2, \mu_L$ and $\mu_R$ are same as in Example \ref{ch2exmp1}. 
 \end{example}

We ran the code for $1000$ iterations and the errors obtained for the final time $\mathscr{T} = 1$ are presented in Table \ref{table:4}. While both PINNs and PA-PINNs yielded poor results for this problem, PA-PINNs performed slightly better in comparison.

\begin{table}[!h]
\centering
\caption{Error comparison for various parameter values for Example \ref{ch2exmp3}.}
\label{table:4}
\begin{tabular}{|c|cc|c|c|c}
\hline \text { Method } & $\varepsilon_{1}$ & $\varepsilon_2$ & $\mathscr{E}_{2}$\\
\hline & $10^{-1}$ & $10^{-2}$ & \textcolor{black}{8.023e-01}\\
\textcolor{black}{\text { PINN }} & $10^{-2}$ & $10^{-3}$ & \textcolor{black}{1.075e+00} \\
& $10^{-3}$ & $10^{-4}$ & \textcolor{black}{8.255e-01}\\
& $10^{-4}$ & $10^{-5}$ & \textcolor{black}{8.777e-01}\\
\hline & $10^{-1}$ & $10^{-2}$ & 1.897e-01 \\
\text { PA-PINN }& $10^{-2}$ & $10^{-3}$ & 5.375e-01 \\
 & $10^{-3}$ & $10^{-4}$ & 6.269e-01 \\
& $10^{-4}$ & $10^{-5}$ & 6.477e-01 \\
\hline
\end{tabular}
\end{table}

\begin{figure}[!h]
  \centering
  \subfloat[Exact Solution vs Predicted solution for ($\varepsilon_1,\varepsilon_2) = (10^{-2},10^{-3})$]{\includegraphics[width=.80\linewidth]{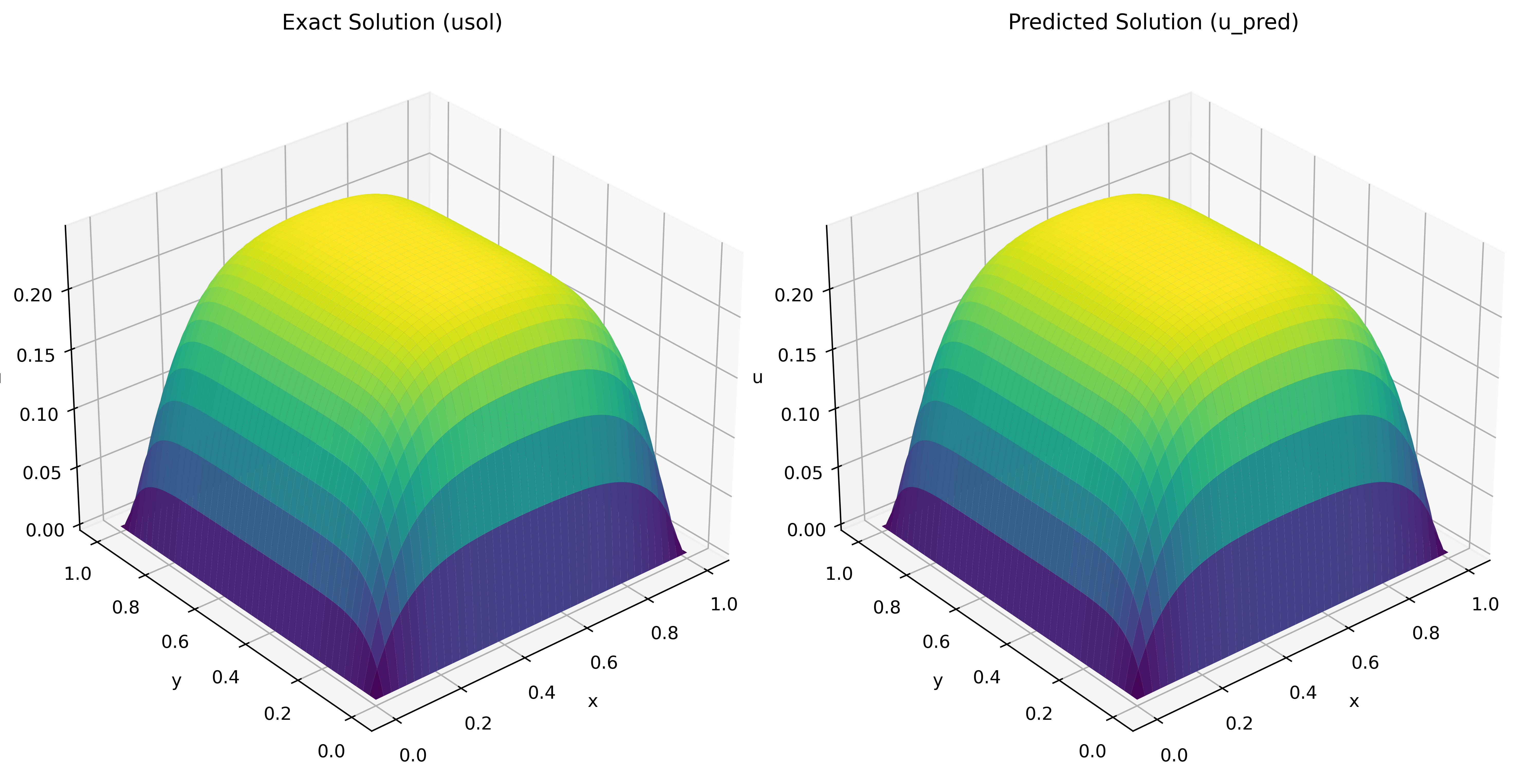}}\\
  \subfloat[Exact Solution vs Predicted Solution for ($\varepsilon_1,\varepsilon_2) = (10^{-3},10^{-4})$]{\includegraphics[width=.80\linewidth]{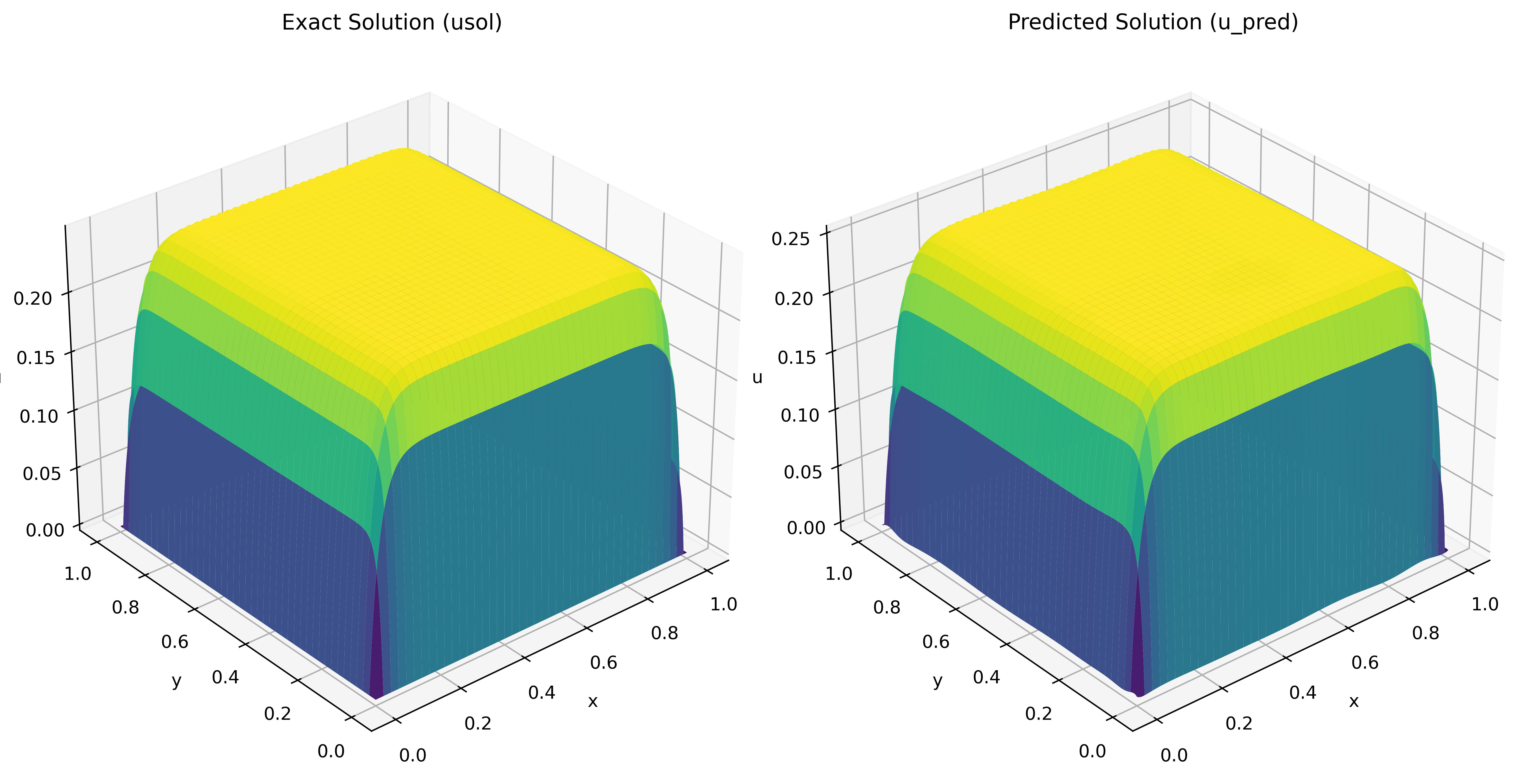}}\\
  \subfloat[Exact Solution vs Predicted Solution for ($\varepsilon_1,\varepsilon_2) = (10^{-4},10^{-5})$]{\includegraphics[width=.80\linewidth]{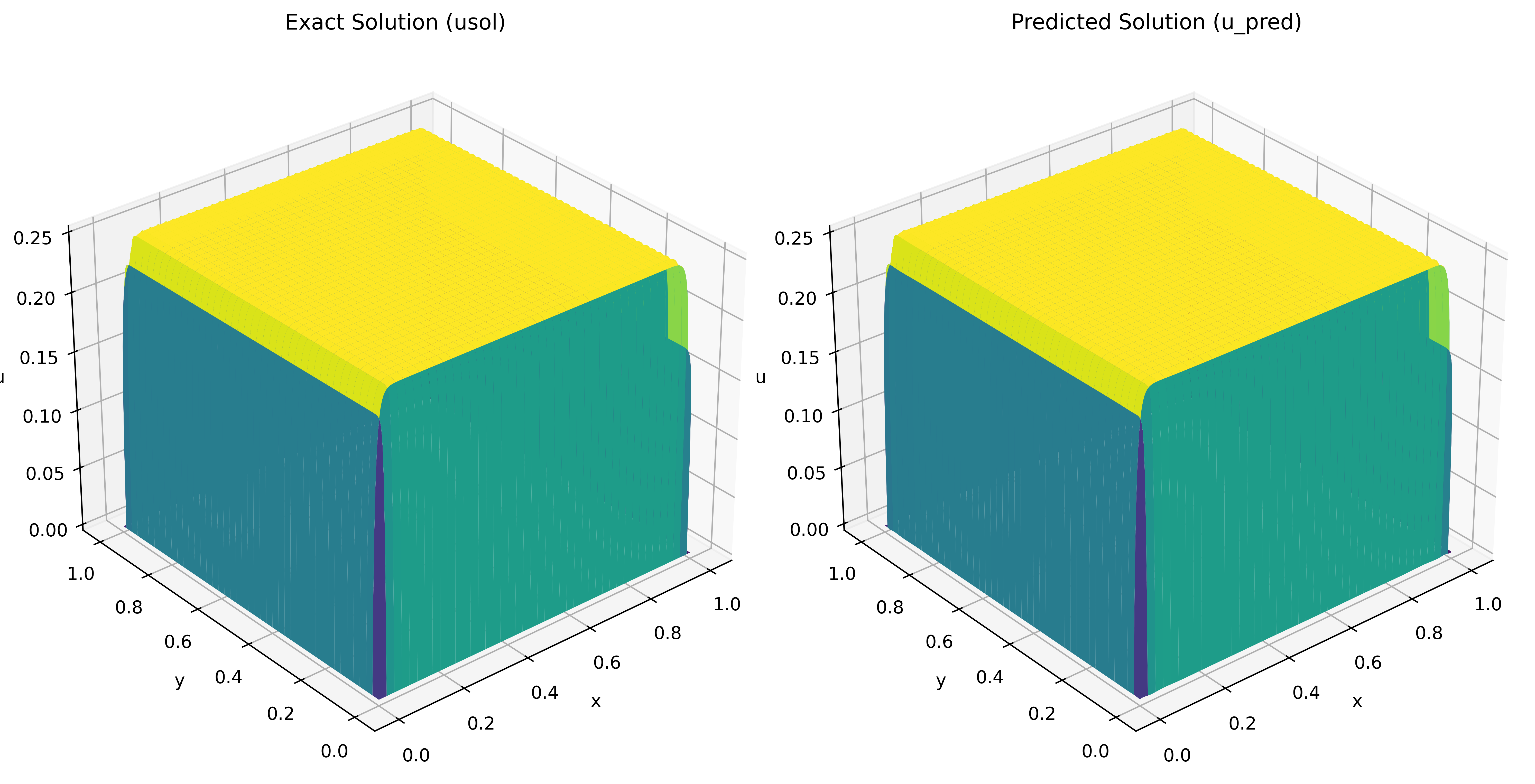}}
   \caption{Comparison between Exact and PA-PINNs solution for different values of perturbation parameter $\varepsilon_1, \varepsilon_2$ for Example \ref{ex2D}.}\label{fig1:2D}
\end{figure}

\begin{example}\label{ex2D}
    Consider the following two-parameter elliptic SPP in 2D
    \begin{equation}\label{ch2_ex4}
		\left\{
		\begin{array}{ll}
		-\varepsilon_1\Delta \mathfrak{v} +\varepsilon_2 \textbf{b}\cdot \nabla \mathfrak{v} +c\mathfrak{v}= g,
		\quad \varOmega\, =  (0,1)^{2},  &\\[6pt]
		\qquad \mathfrak{v} = 0\, , \quad (x,y) \in \partial \varOmega,
		\end{array} \right.
		\end{equation}
  where \textbf{b} $= (b_1(x,y), b_2(x,y)) = (1,0)$ and $c = 1$.
  We choose $g$ in a manner that
  \begin{eqnarray*}
  \mathfrak{v}(x,y) &=& 0.25\times\Big(1-\exp(-\varepsilon_2 l_{1} x/2\varepsilon_1)\Big) \times \Big(1-\exp(-y/\sqrt{\varepsilon_1})\Big)\\
  && \qquad \times \Big(1-\exp(-\varepsilon_2 {l}_{2} (1-x)/2\varepsilon_1)\Big)\times \Big(1-\exp(-(1-y)/\sqrt{\varepsilon_1})\Big),
  \end{eqnarray*}
  is the exact solution and
  \[
  l_{1,2} = \left(\sqrt{1+16\dfrac{\varepsilon_1}{\varepsilon_2^{2}}}\right) \mp 1.
  \]
\end{example}
The solution of the elliptic BVP \eqref{ch2_ex4} along $x = 0$ and $x = 1$ has exponential boundary layers. Also, along $y = 0$ and $y = 1$ characteristic boundary layers can be observed along with the corner layers at each corner of domain $\varOmega$. Figure \ref{fig1:2D} shows that PA-PINN approximates the solution with sharp boundary layers quite easily. Table \ref{table:5} contains the relative $L^{2}-$errors obtained for various values of perturbation parameters.

\begin{table}[!h]
\centering
\caption{Error comparison for various parameter values for Example \ref{ex2D}.}
\label{table:5}
\begin{tabular}{|c|cc|c|c|c}
\hline \text { Method } & $\varepsilon_{1}$ & $\varepsilon_2$ & $\mathscr{E}_{2}$ \\
\hline & $10^{-1}$ & $10^{-2}$ & \textcolor{black}{4.418e-02}\\
\textcolor{black}{\text { PINN }} & $10^{-2}$ & $10^{-3}$ & \textcolor{black}{8.363e-02} \\
& $10^{-3}$ & $10^{-4}$ & \textcolor{black}{1.520e-01}\\
& $10^{-4}$ & $10^{-5}$ & \textcolor{black}{2.467e-01}\\
\hline & $10^{-1}$ & $10^{-2}$ & 1.011e-04 \\
\text { PA-PINN }& $10^{-2}$ & $10^{-3}$ & 1.084e-04 \\
 & $10^{-3}$ & $10^{-4}$ &  2.621e-04  \\
& $10^{-4}$ & $10^{-5}$ & 3.812e-04  \\
\hline
\end{tabular}
\end{table}

\begin{figure}[!h]
  \centering
  \subfloat[Exact Solution vs Predicted Solution for ($\varepsilon_1,\varepsilon_2) = (10^{-1},10^{-2})$]{\includegraphics[width=.80\linewidth]{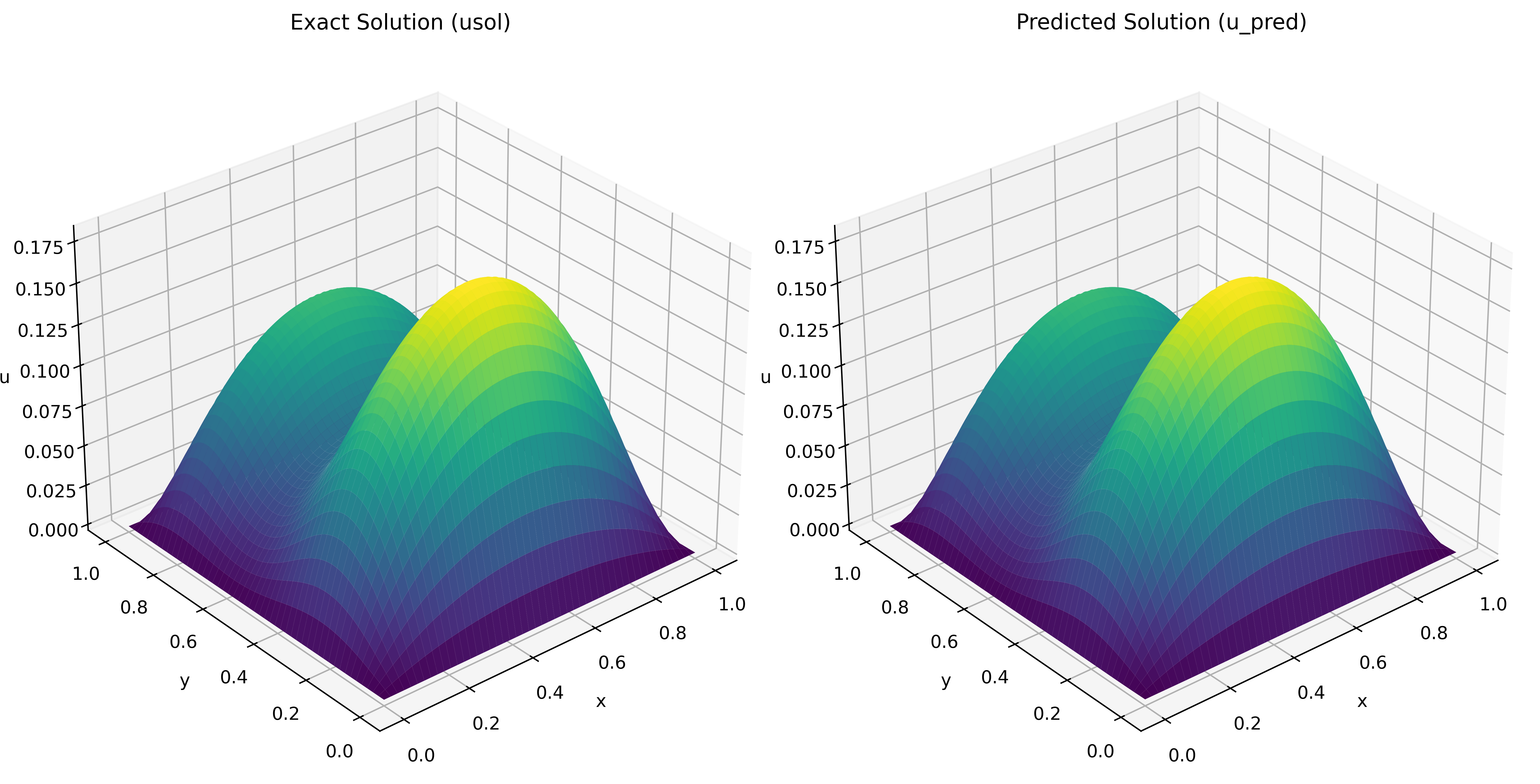}}\\
    \subfloat[Exact Solution vs Predicted Solution for ($\varepsilon_1,\varepsilon_2) = (10^{-2},10^{-3})$]{\includegraphics[width=.80\linewidth]{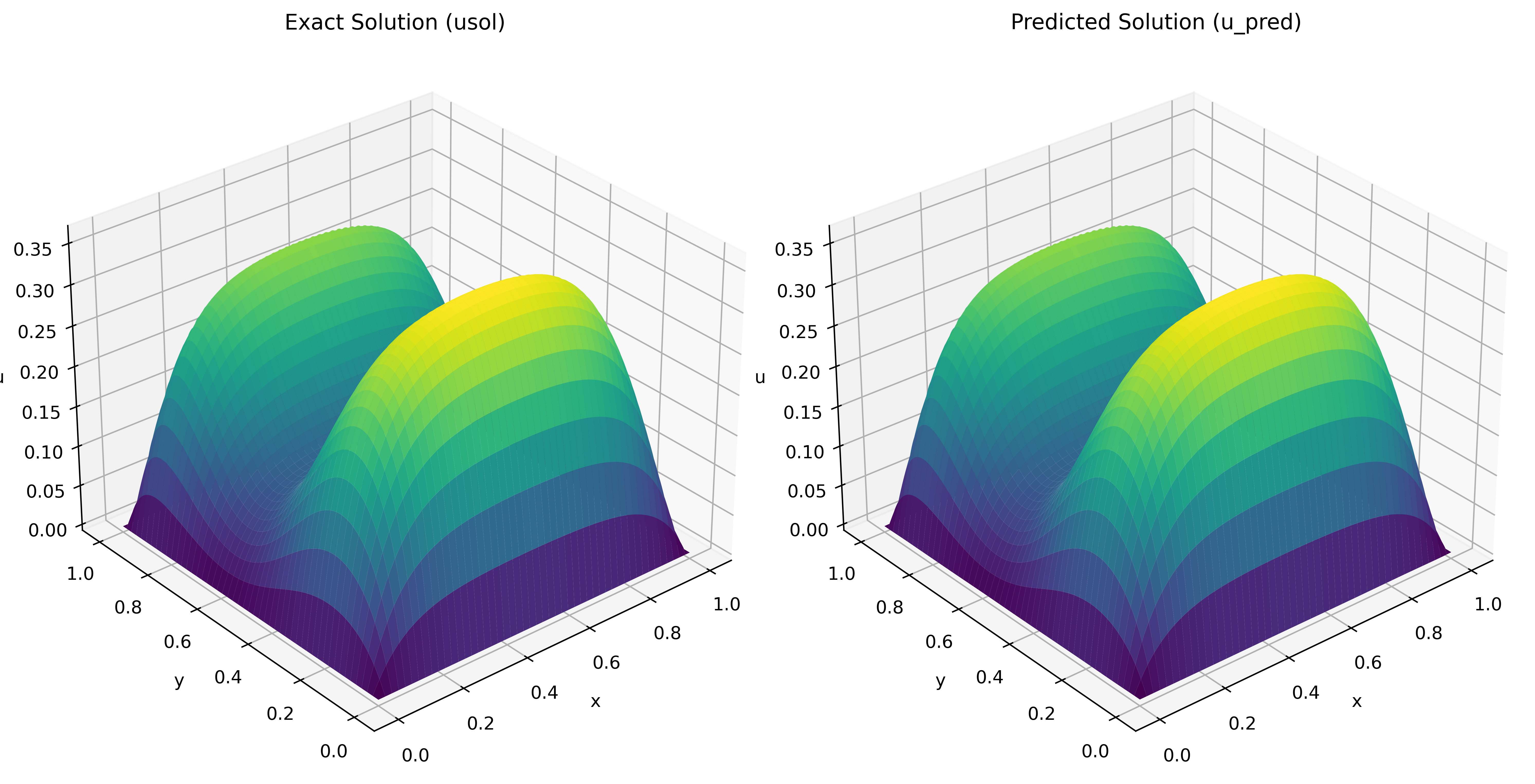}}\\
     \subfloat[Exact Solution vs Predicted Solution for ($\varepsilon_1,\varepsilon_2) = (10^{-3},10^{-4})$]{\includegraphics[width=.80\linewidth]{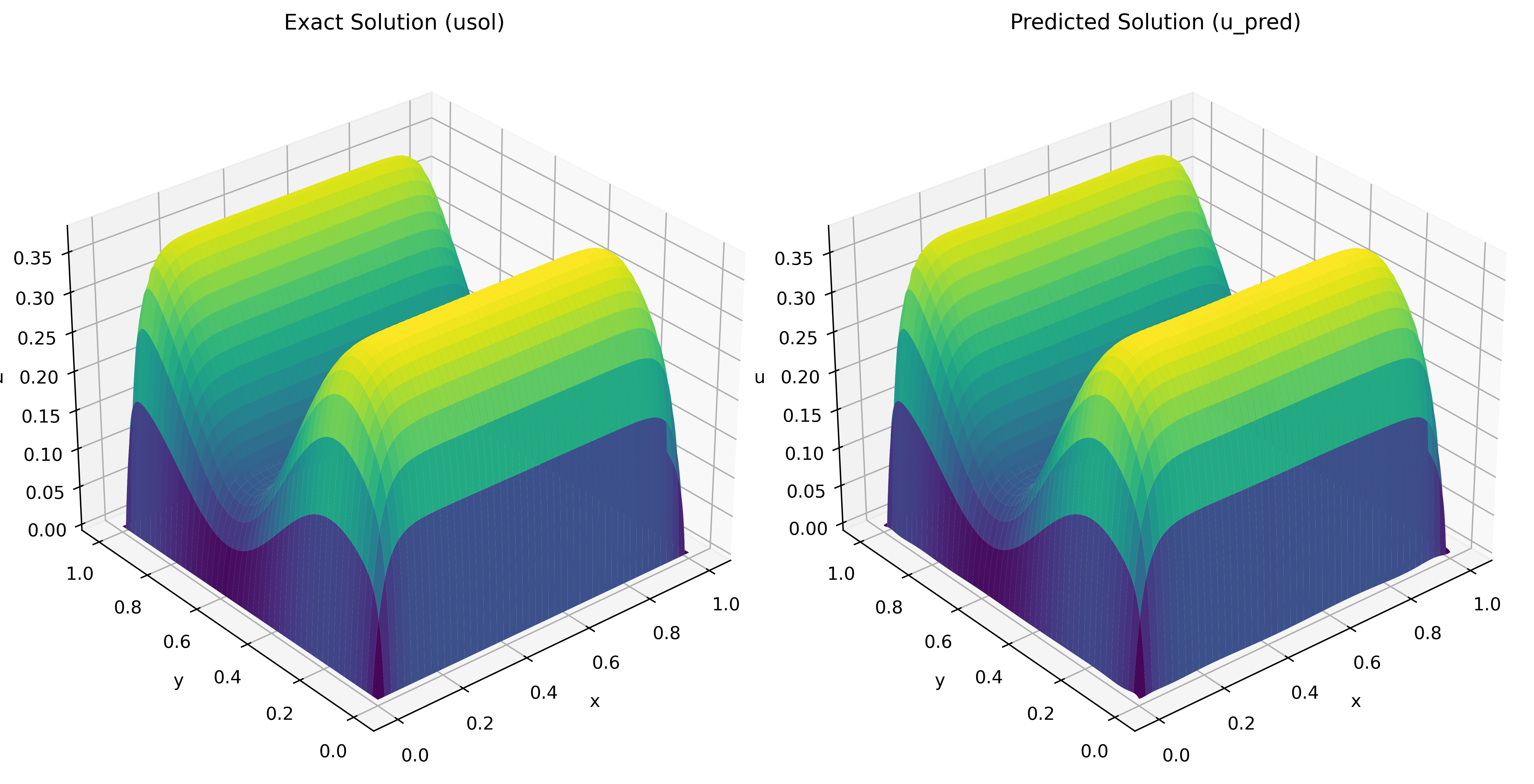}}
  \caption{Comparison between Exact and PA-PINNs solution for different values of perturbation parameter $\varepsilon_1, \varepsilon_2$ for Example \ref{ex2D_2}.}\label{fig2:2D}
\end{figure}

\begin{example}\label{ex2D_2}
    Consider the following test problem
    \begin{equation}\label{}
		\left\{
		\begin{array}{ll}
		-\varepsilon_1\Delta \mathfrak{v} +\varepsilon_2 \mathfrak{v}_{x} + \mathfrak{v} = h(x,y),
		\quad \varOmega\, =  (0,1)^{2},  &\\[6pt]
		\qquad \mathfrak{v} = 0\, , \quad (x,y) \in \partial \varOmega .
		\end{array} \right.
		\end{equation}
 Choosing $h(x,y)$ in such a way that
  \begin{equation*}
  \mathfrak{v}(x,y) = \dfrac{1}{4}\left(1+\dfrac{\sin(8x)}{2}\right) \times \mathfrak{w}(x,y),
  \end{equation*}
  is the exact solution where
  \begin{eqnarray*}
      \mathfrak{w}(x,y) &=& \Big(1-\exp(-\varepsilon_2 l_{1} x/2\varepsilon_1)\Big) \times \Big(1-\exp(-y/\sqrt{\varepsilon_1})\Big)\\
      &&\times \Big(1-\exp(-\varepsilon_2 {l}_{2} (1-x)/2\varepsilon_1)\Big) \times \Big(1-\exp(-(1-y)/\sqrt{\varepsilon_1})\Big),
  \end{eqnarray*}
 and $l_{1,2}$ is same as in Example \ref{ex2D}.
\end{example}
 Error values are presented in Table \ref{table:6}. Figure \ref{fig2:2D} shows the exact and PA-PINNs solution plot comparison for different values of perturbation parameters $\varepsilon_1$ and $\varepsilon_2$. It is evident that the PA-PINN solution aligns well with the exact solution across various parameter values.   

\begin{table}[ht]
\centering
\caption{Error comparison for various parameter values for Example \ref{ex2D_2}.}
\label{table:6}
\begin{tabular}{|c|cc|c|c|c}
\hline \text { Method } & $\varepsilon_{1}$ & $\varepsilon_2$ & $\mathscr{E}_{2}$ \\
\hline & $10^{-1}$ & $10^{-2}$ & \textcolor{black}{4.680e-02}\\
\textcolor{black}{\text { PINN }} & $10^{-2}$ & $10^{-3}$ & \textcolor{black}{9.387e-02} \\
& $10^{-3}$ & $10^{-4}$ & \textcolor{black}{2.352e-01}\\
& $10^{-4}$ & $10^{-5}$ & \textcolor{black}{2.962e-01}\\
\hline & $10^{-1}$ & $10^{-2}$ & 5.609e-04 \\
\text { PA-PINN }& $10^{-2}$ & $10^{-3}$ & 1.895e-04 \\
 & $10^{-3}$ & $10^{-4}$ & 2.232e-03 \\
& $10^{-4}$ & $10^{-5}$ & 8.933e-03 \\
\hline
\end{tabular}
\end{table}

\begin{example}\label{ex2D_time}
    Here we will be considering parabolic two-parameter problem which is  similar to Example \ref{ex2D}.
    \begin{equation}\label{ch2_ex5}
		\left\{
		\begin{array}{ll}
		\dfrac{\partial \mathfrak{v}}{\partial t}-\varepsilon_1\Delta \mathfrak{v} +\varepsilon_2 \textbf{b}\cdot \nabla \mathfrak{v} +c\mathfrak{v}= g,
		\quad \varOmega\, =  (0,1)^{2},\, 0\leq t\leq \mathscr{T}, &\\ [6pt]
		\mathfrak{v}(x,y,0) = \mathfrak{v}_{0}(x,y),\quad \mathfrak{v}(x,y,t) = 0\, , \quad (x,y) \in \partial \varOmega,
		\end{array} \right.
		\end{equation}
  where \textbf{b} $= (b_1(x,y), b_2(x,y)) = (1,0)$ and $c = 1$.
  We choose $g$ in a manner that
  \begin{eqnarray*}
  \mathfrak{v}(x,y) &=& 0.25\times\Big(1-\exp(-t)\Big)\times\Big(1-\exp(-\varepsilon_2 l_{1} x/2\varepsilon_1)\Big) \times \Big(1-\exp(-y/\sqrt{\varepsilon_1})\Big)\\
  && \qquad\times\, \Big(1-\exp(-\varepsilon_2 {l}_{2} (1-x)/2\varepsilon_1)\Big)\times \Big(1-\exp(-(1-y)/\sqrt{\varepsilon_1})\Big),
  \end{eqnarray*}
  is the exact solution and
  \[
  l_{1,2} = \left(\sqrt{1+16\dfrac{\varepsilon_1}{\varepsilon_2^{2}}}\right) \mp 1.
  \]
\end{example}
$\mathscr{E}_{2}$ and $\mathscr{E}_{\infty}$ error values for the final time $\mathscr{T} = 1$ are  presented in Table \ref{table:7} for various values of perturbation parameters $\varepsilon_1$ and $\varepsilon_2$. PA-PINN once again outperforms PINNs in terms of accuracy, as observed in previous examples as well.

\begin{table}[ht]
\centering
\caption{Error comparison for various parameter values for Example \ref{ex2D_time}.}
\label{table:7}
\begin{tabular}{|c|cc|c|c|c}
\hline \text { Method } & $\varepsilon_{1}$ & $\varepsilon_2$ & $\mathscr{E}_{2}$ & $\mathscr{E}_{\infty}$\\
\hline & $10^{-1}$ & $10^{-2}$ & \textcolor{black}{4.944e-01} & \textcolor{black}{2.999e-02}\\
\textcolor{black}{\text {PINN }}& $10^{-2}$ & $10^{-3}$ & \textcolor{black}{4.492e-01} & \textcolor{black}{6.970e-02}\\
 & $10^{-3}$ & $10^{-4}$ & \textcolor{black}{4.730e-01} & \textcolor{black}{9.653e-02}\\
& $10^{-4}$ & $10^{-5}$ & \textcolor{black}{5.465e-01} & \textcolor{black}{1.192e-01}\\
\hline & $10^{-1}$ & $10^{-2}$ & 6.941e-04 & 2.363e-04\\
\text { PA-PINN }& $10^{-2}$ & $10^{-3}$ & 4.626e-04 & 2.918e-04\\
 & $10^{-3}$ & $10^{-4}$ & 1.417e-03 & 1.130e-03\\
& $10^{-4}$ & $10^{-5}$ & 1.173e-03 & 4.361e-03\\
\hline
\end{tabular}
\end{table}

\subsection{Accuracy of PA-PINNs with different \textcolor{black}{width, depth and training points}}
\begin{figure}[!h]
  \centering
  \subfloat[With 1000 training points]{\includegraphics[width=.33\linewidth]{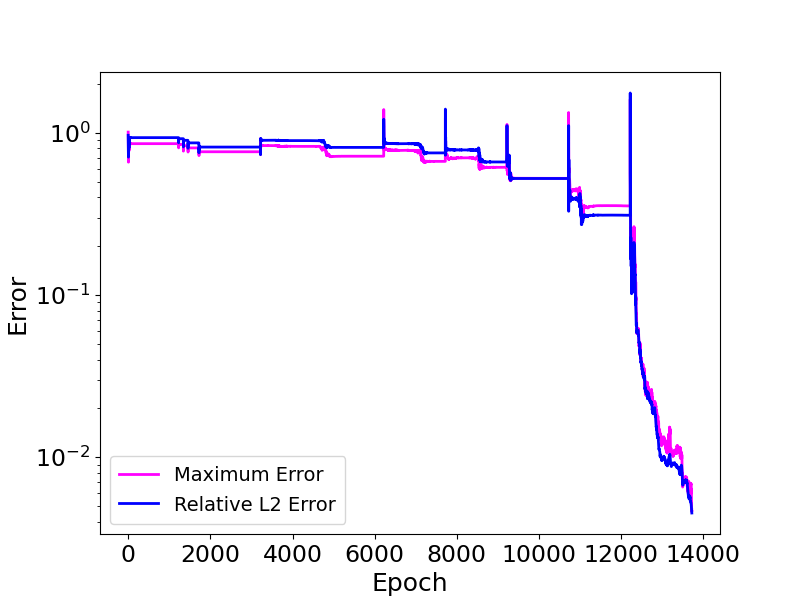}}\hfill
  \subfloat[With 5000 training points] {\includegraphics[width=.33\linewidth]{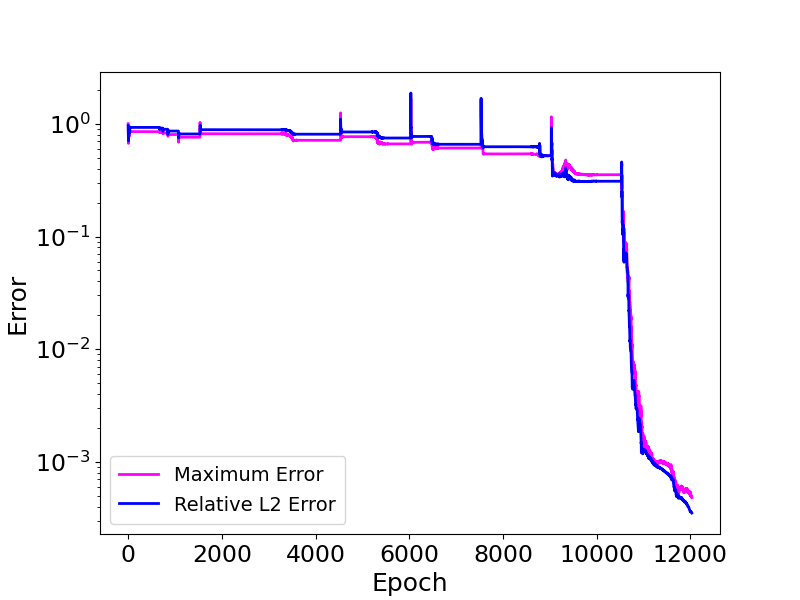}}
  \subfloat[With 10000 training points] {\includegraphics[width=.33\linewidth]{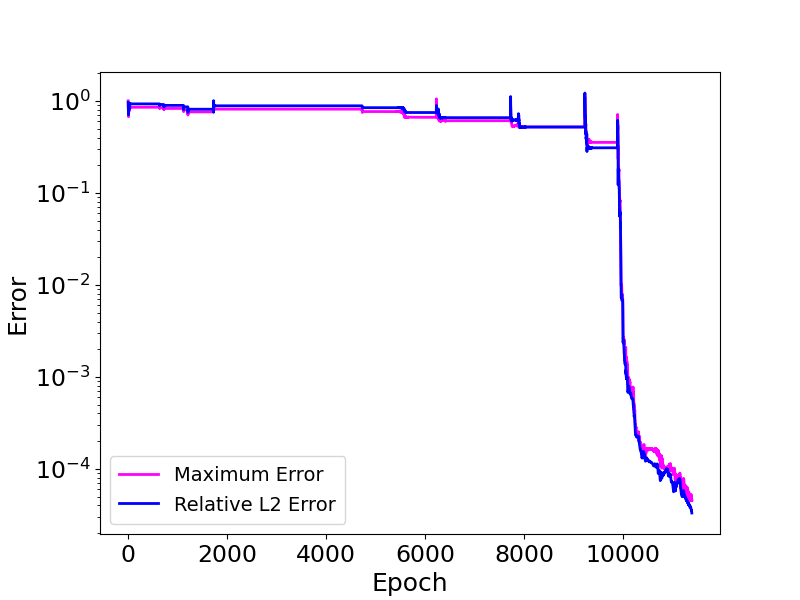}}
  \caption{Error vs epoch plots for Example \ref{ch2exmp1} with $(\varepsilon_1, \varepsilon_2) = (10^{-3}, 10^{-4})$.}\label{fig11_ex1}
\end{figure}

\begin{figure}[!h]
  \centering
  \subfloat[With 1000 training points]{\includegraphics[width=.33\linewidth]{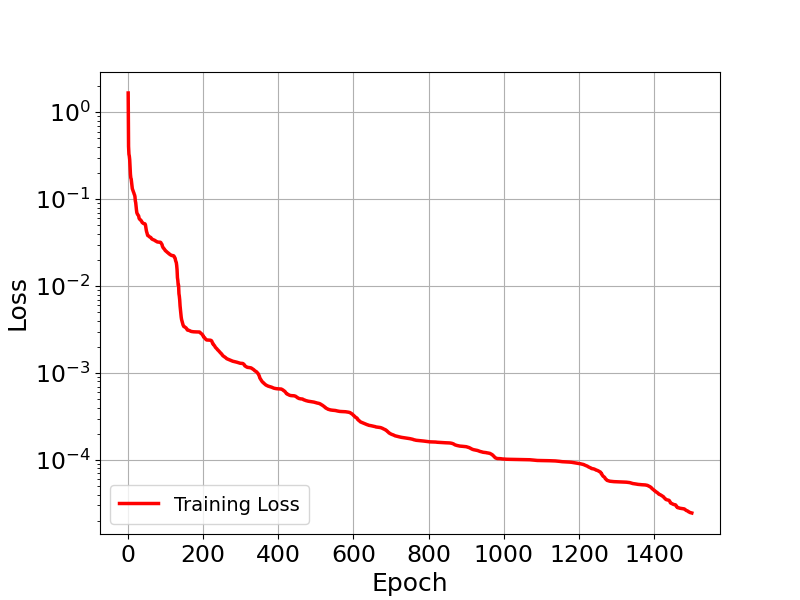}}\hfill
  \subfloat[With 5000 training points] {\includegraphics[width=.33\linewidth]{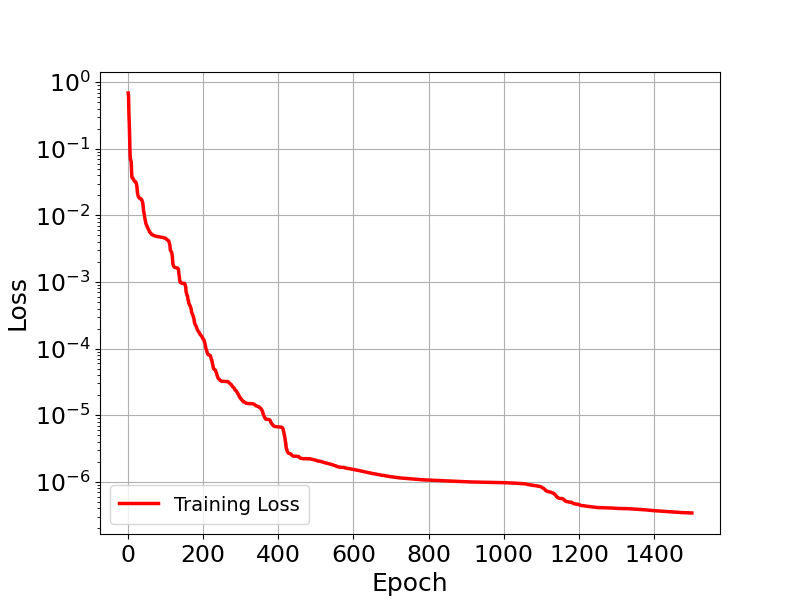}}
  \subfloat[With 10000 training points] {\includegraphics[width=.33\linewidth]{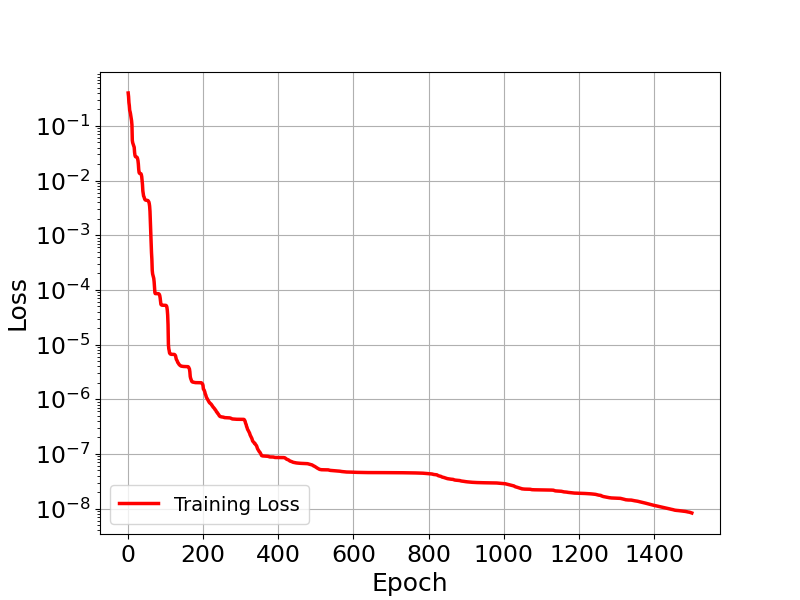}}
  \caption{Loss curve for Example \ref{ch2exmp1} with $(\varepsilon_1, \varepsilon_2) = (10^{-3}, 10^{-4})$.}\label{fig111_ex1}
\end{figure}

\begin{figure}[!h]
  \centering
  \subfloat[With 1000 training points]{\includegraphics[width=.33\linewidth]{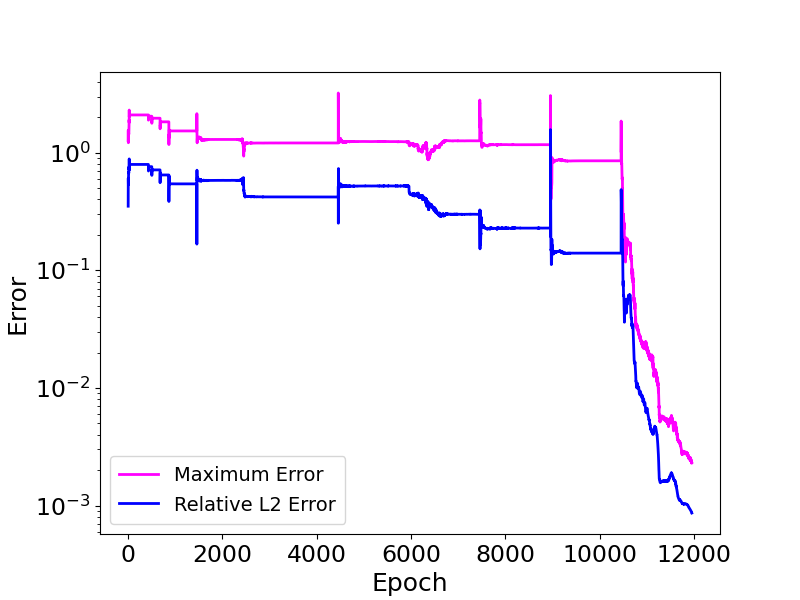}}\hfill
  \subfloat[With 5000 training points] {\includegraphics[width=.33\linewidth]{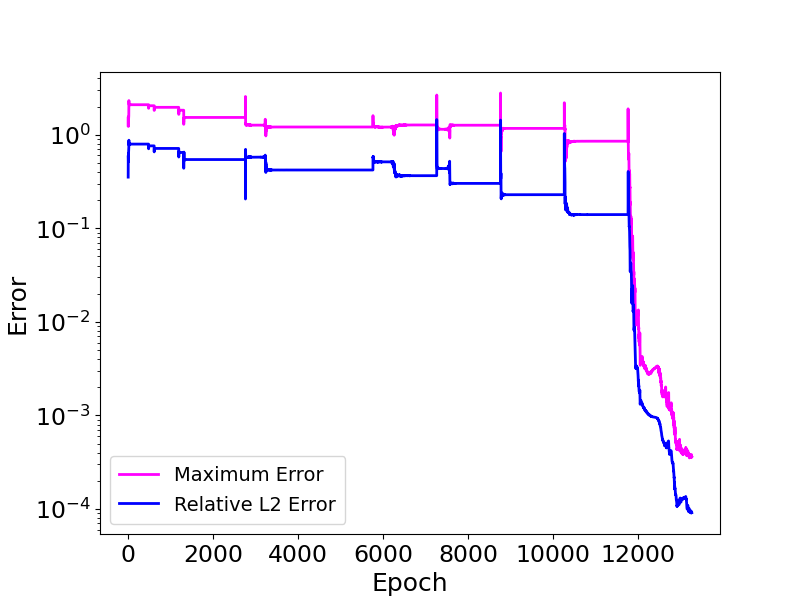}}
  \subfloat[With 10000 training points] {\includegraphics[width=.33\linewidth]{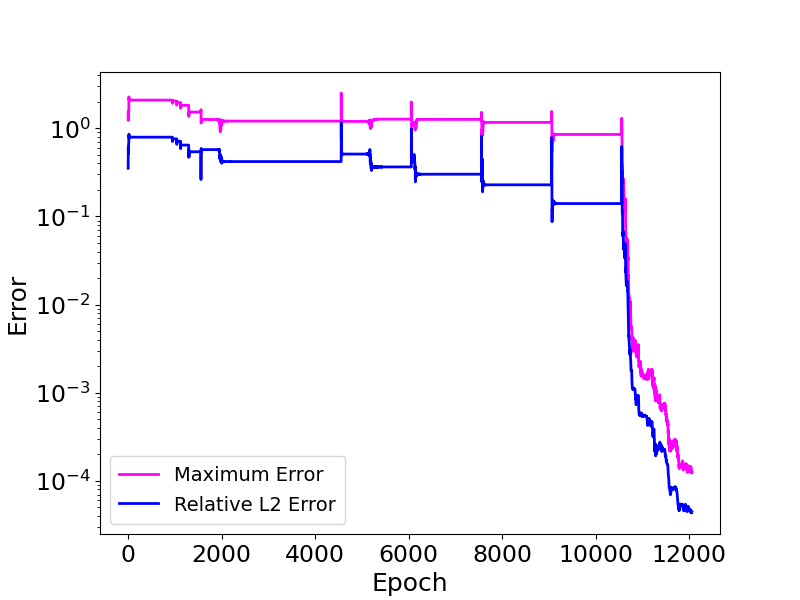}}
  \caption{Error vs epoch plots for Example \ref{ch2exmp2} with $(\varepsilon_1, \varepsilon_2) = (10^{-3}, 10^{-4})$.}\label{err_epoch_ex2}
\end{figure}

\begin{figure}[!h]
  \centering
  \subfloat[With 1000 training points]{\includegraphics[width=.33\linewidth]{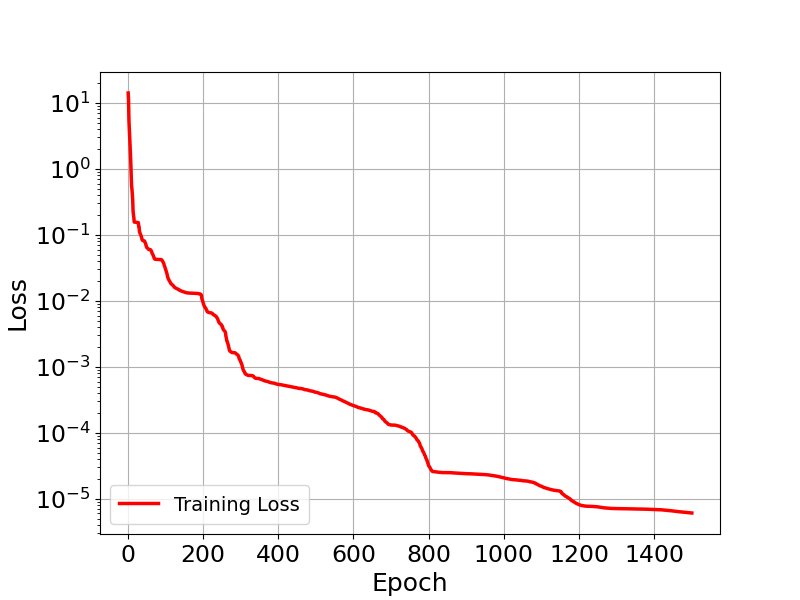}}\hfill
  \subfloat[With 5000 training points] {\includegraphics[width=.33\linewidth]{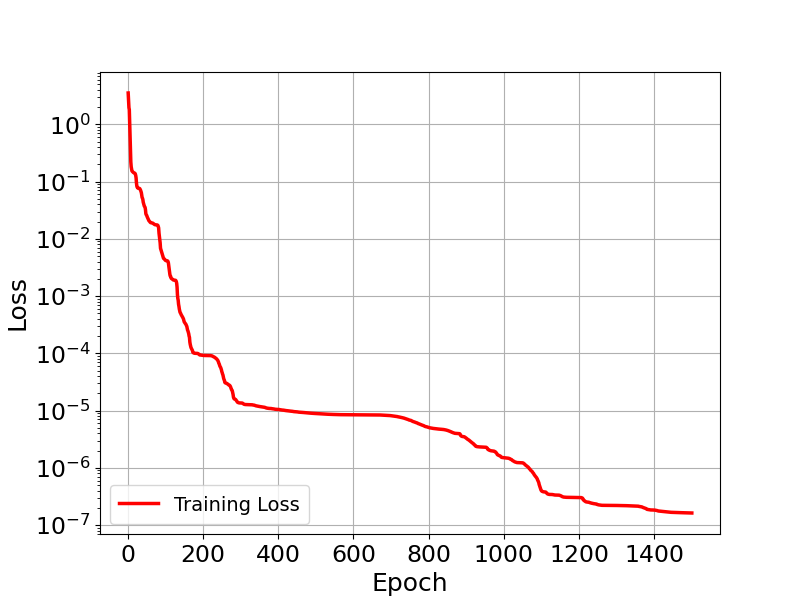}}
  \subfloat[With 10000 training points] {\includegraphics[width=.33\linewidth]{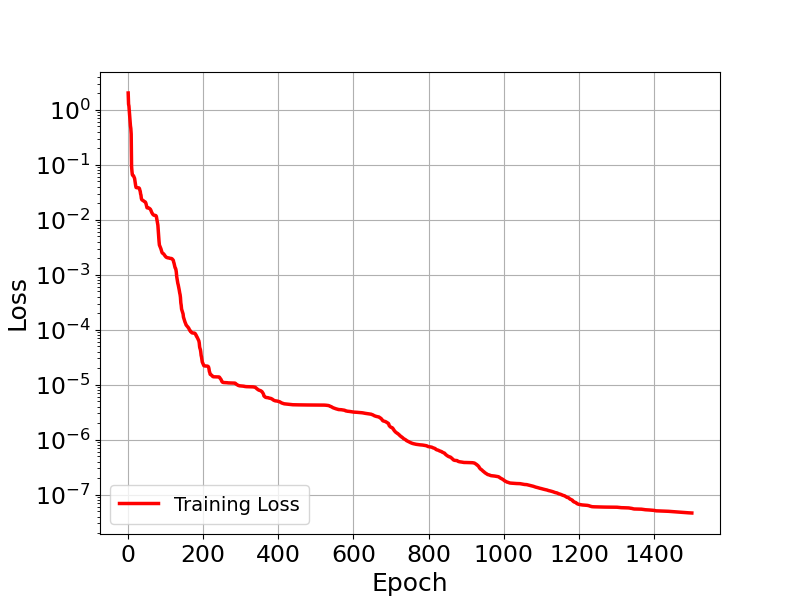}}
  \caption{Loss curve for Example \ref{ch2exmp2} with $(\varepsilon_1, \varepsilon_2) = (10^{-3}, 10^{-4})$.}\label{loss_ex2}
\end{figure}

\begin{figure}[!h]
  \centering
  \subfloat[With 1000 training points]{\includegraphics[width=.33\linewidth]{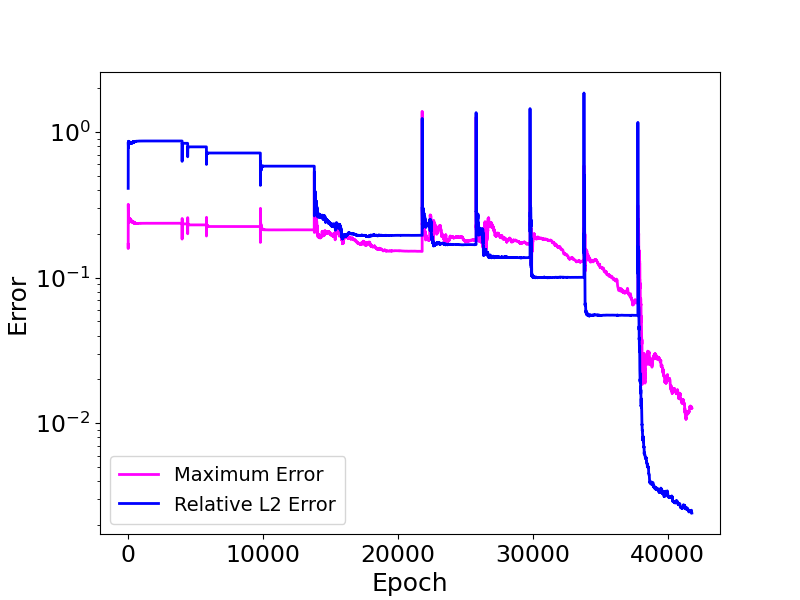}}\hfill
  \subfloat[With 5000 training points] {\includegraphics[width=.33\linewidth]{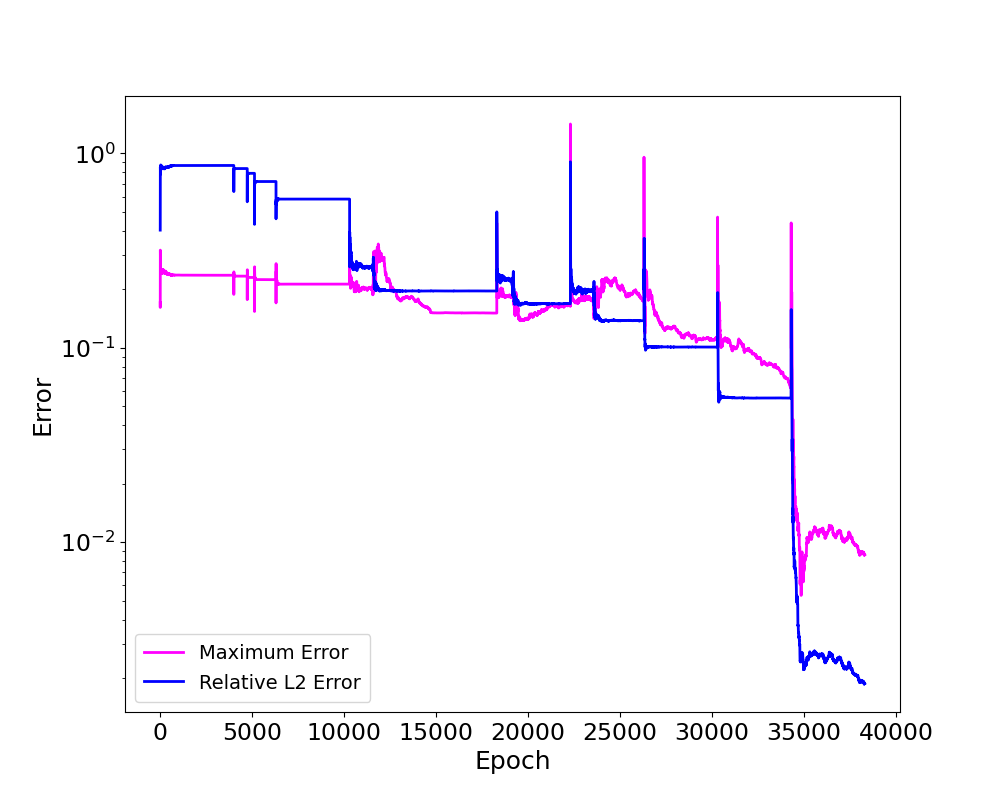}}
  \subfloat[With 10000 training points] {\includegraphics[width=.33\linewidth]{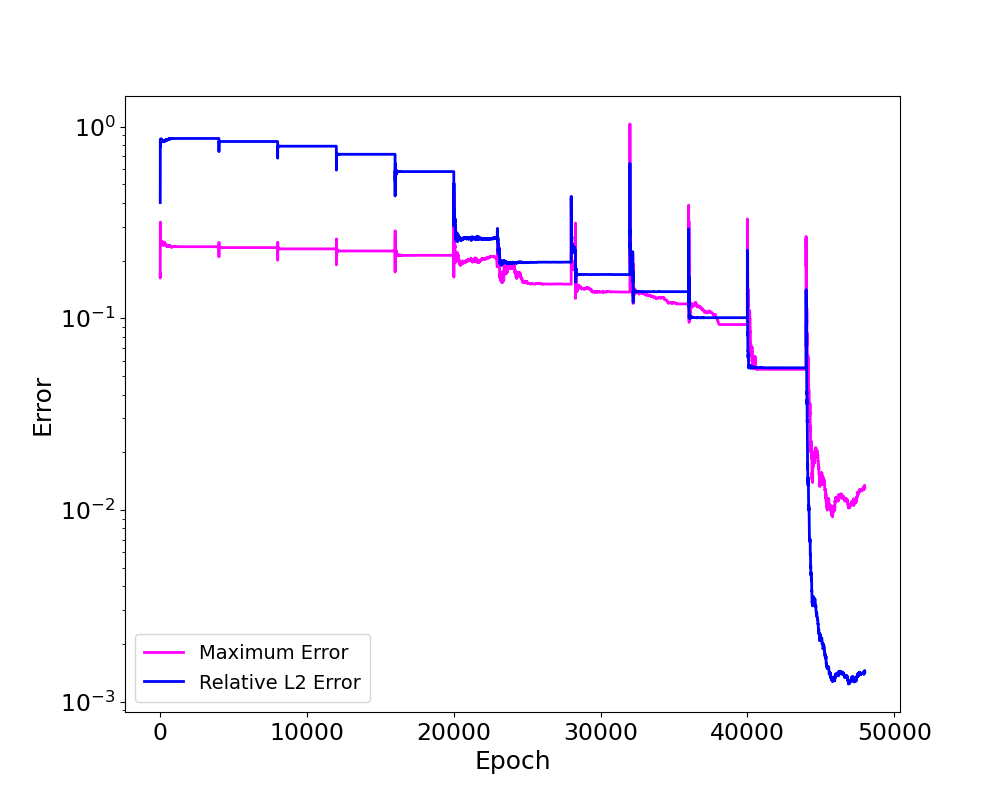}}
  \caption{Error vs epoch plots for Example \ref{ex2D} with $(\varepsilon_1, \varepsilon_2) = (10^{-3}, 10^{-4})$.}\label{err_epoch_ex2D}
\end{figure}

\begin{figure}[!h]
  \centering
  \subfloat[With 1000 training points]{\includegraphics[width=.33\linewidth]{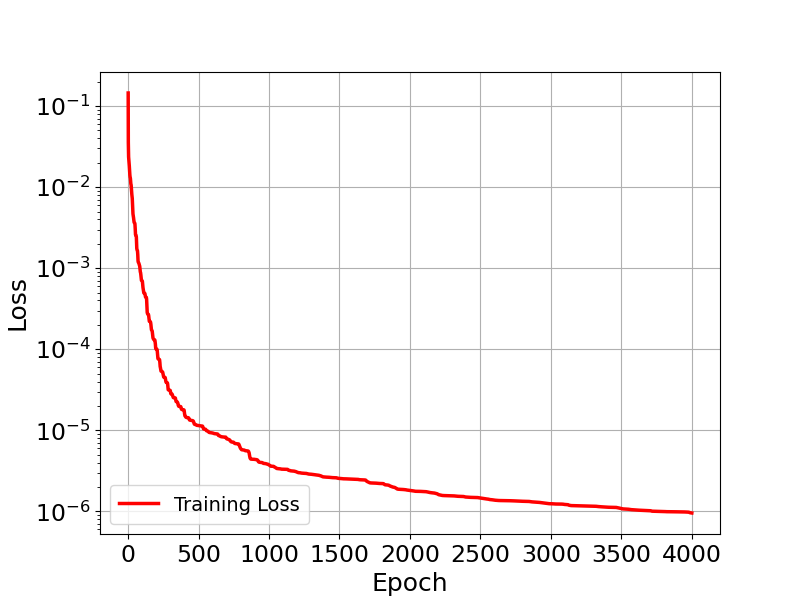}}\hfill
  \subfloat[With 5000 training points] {\includegraphics[width=.33\linewidth]{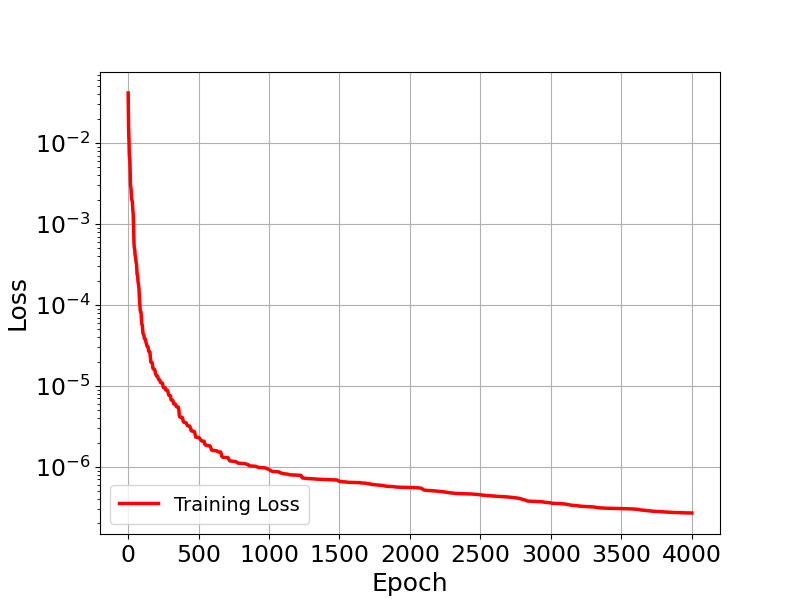}}
  \subfloat[With 10000 training points] {\includegraphics[width=.33\linewidth]{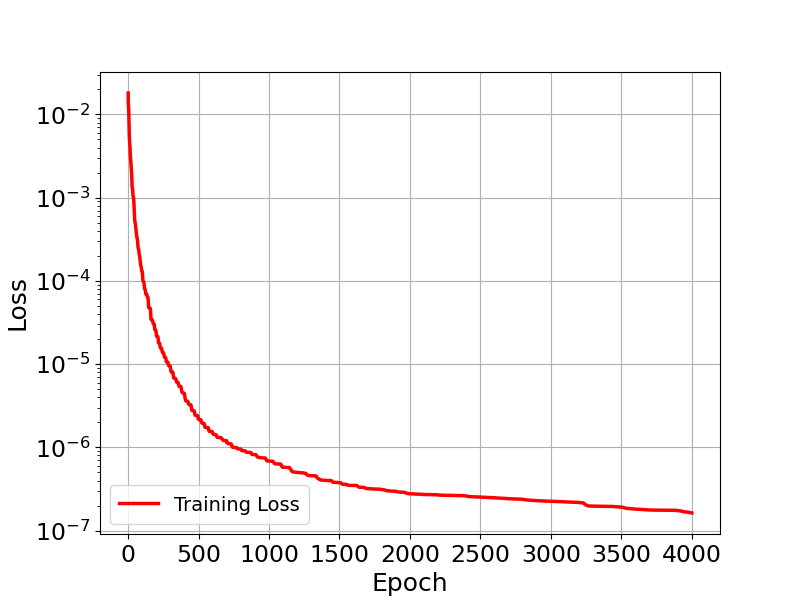}}
  \caption{Loss curve for Example \ref{ex2D} with $(\varepsilon_1, \varepsilon_2) = (10^{-3}, 10^{-4})$.}\label{loss_ex2D}
\end{figure}

\begin{figure}[!h]
  \centering
  \subfloat[With 1000 training points]{\includegraphics[width=.33\linewidth]{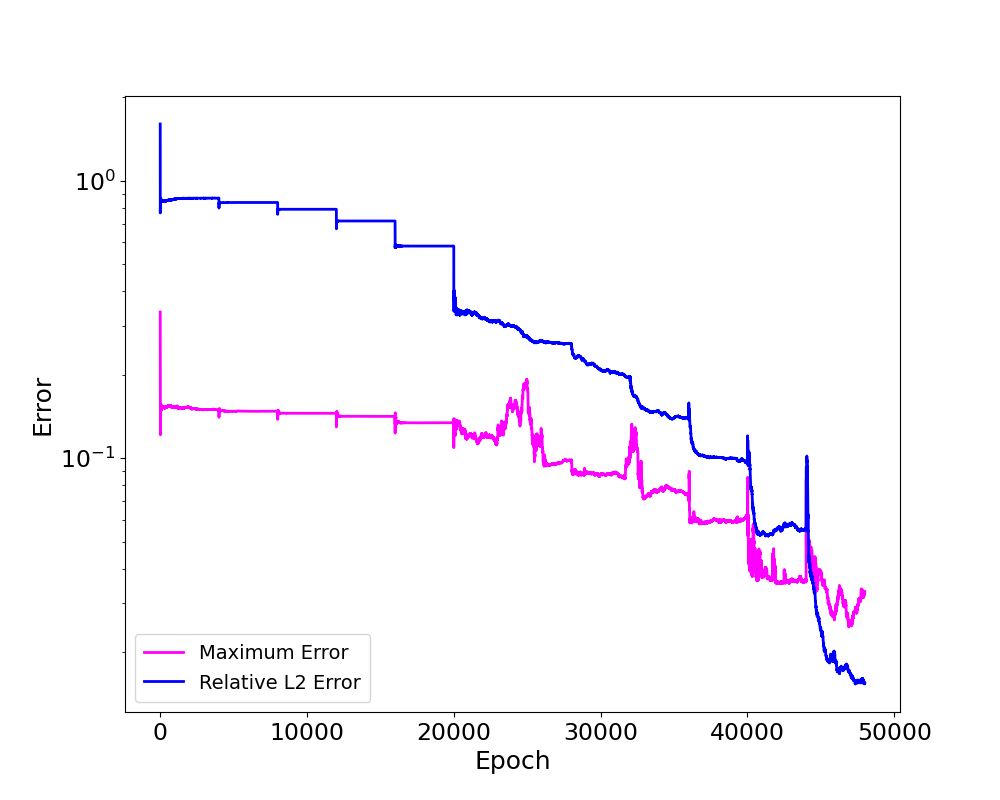}}\hfill
  \subfloat[With 5000 training points] {\includegraphics[width=.33\linewidth]{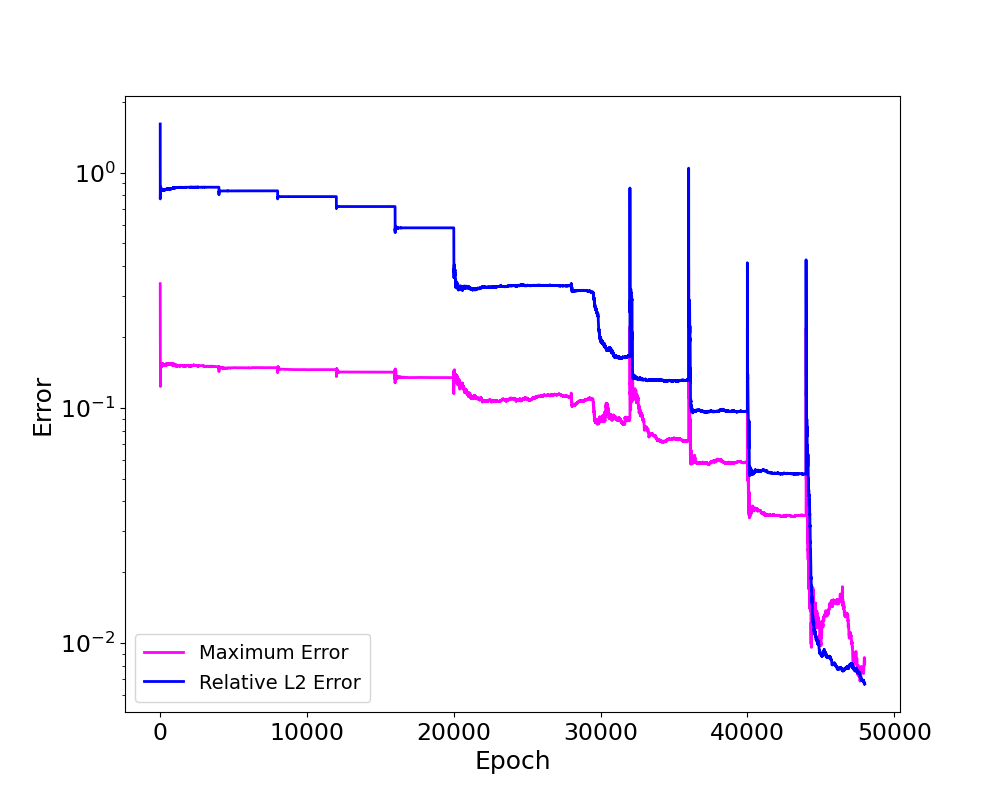}}
  \subfloat[With 10000 training points] {\includegraphics[width=.33\linewidth]{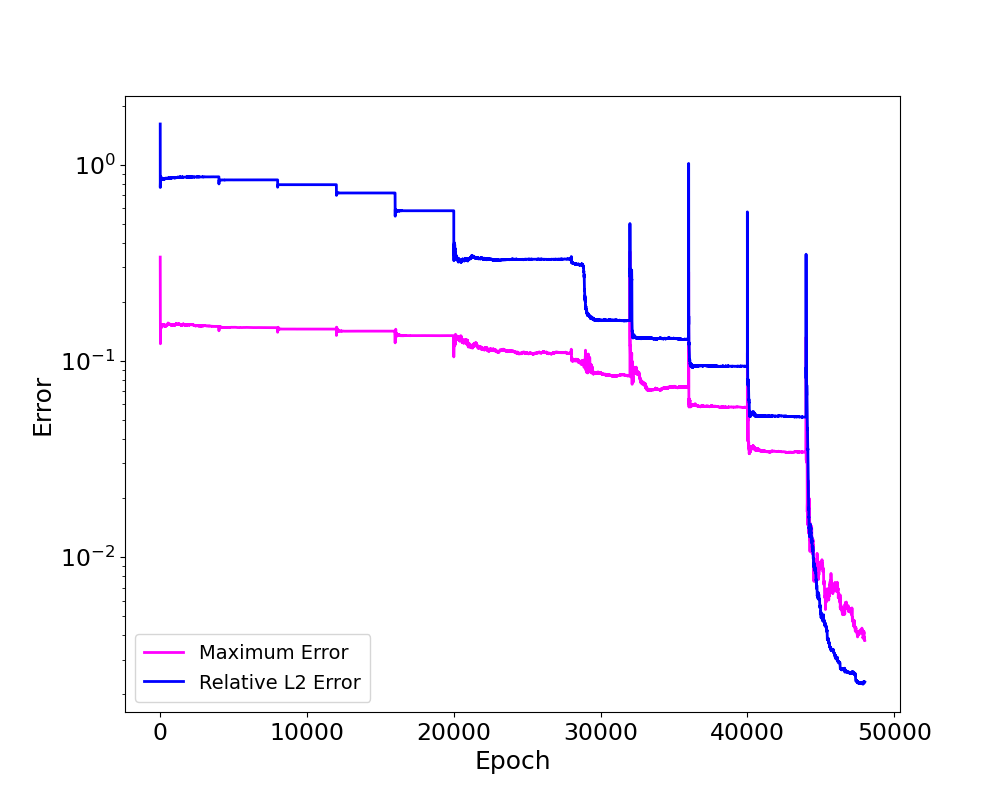}}
  \caption{Error vs epoch plots for Example \ref{ex2D_time} with $(\varepsilon_1, \varepsilon_2) = (10^{-3}, 10^{-4})$.}\label{err_epoch_ex2D_time}
\end{figure}

\begin{figure}[!h]
  \centering
  \subfloat[With 1000 training points]{\includegraphics[width=.33\linewidth]{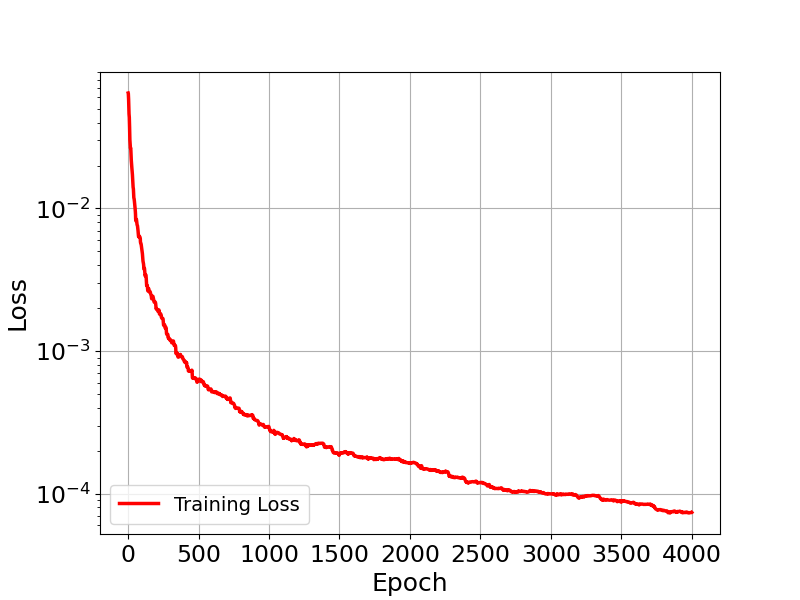}}\hfill
  \subfloat[With 5000 training points] {\includegraphics[width=.33\linewidth]{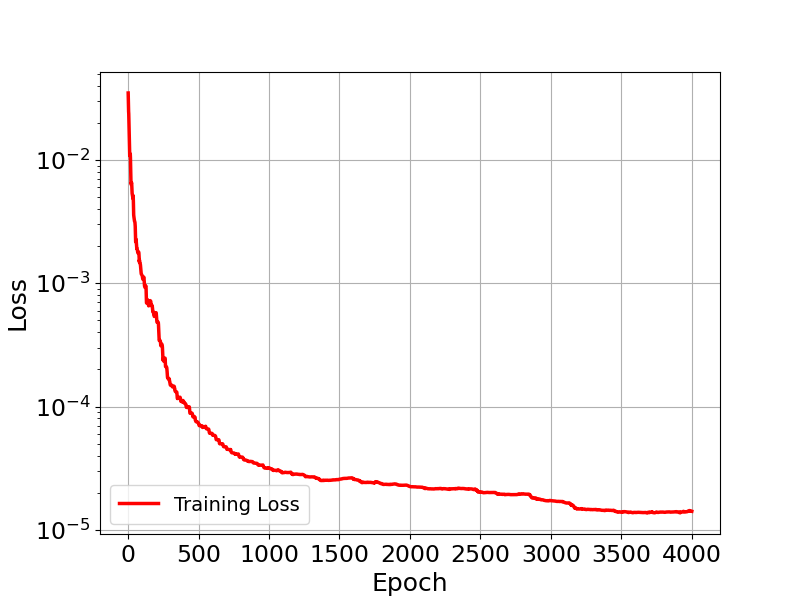}}
  \subfloat[With 10000 training points] {\includegraphics[width=.33\linewidth]{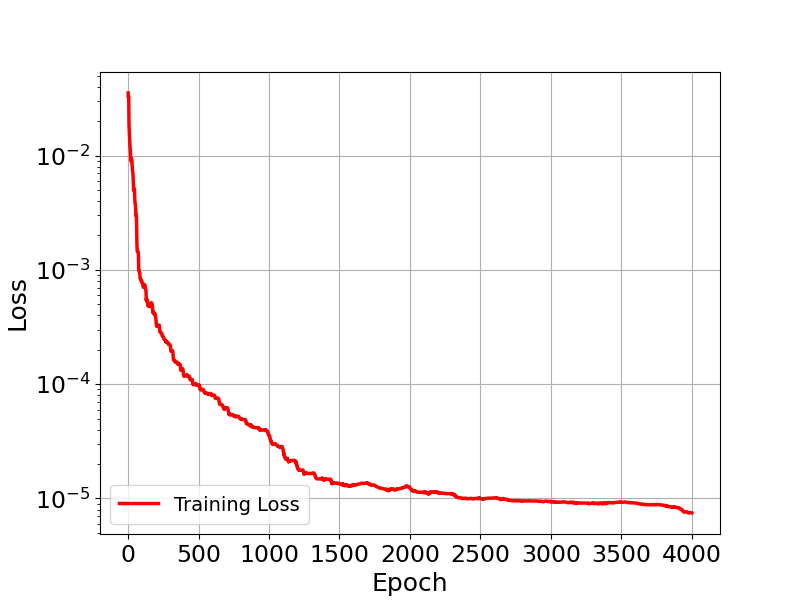}}
  \caption{Loss curve for Example \ref{ex2D_time} with $(\varepsilon_1, \varepsilon_2) = (10^{-3}, 10^{-4})$.}\label{loss_ex2D_time}
\end{figure}

Here, we studied the impact of varying neural network width, depth and training points on the accuracy and computational cost. We have considered network with different depth and width for comparison, rest of the network hyperparameters are as it is. We have taken $(\varepsilon_{1}^{0} , \varepsilon_{2}^{0}) = (0.3 , 0.4)$, parameter values as $c_1 = 0.71$ and $c_2 = 0.72$. Tables \ref{non_linear_1} and \ref{non_linear_2} display the $ L^{2}$ errors for Example \ref{ch2exmp1} and Example \ref{ch2exmp2}, respectively, computed for various widths and depths using  nonlinear parameter update. Error values for Example \ref{ch2exmp1} and Example \ref{ch2exmp2}, obtained using the linear parameter update, are shown in Tables \ref{linear_1} and \ref{linear_2}, respectively. Table \ref{linear_3} and Table \ref{ex2D_table} shows the error obtained by using linear-update rule for Example \ref{ch2exmp3} and \ref{ex2D} respectively. Considering the accuracy and computational cost, we choose size of network to be $8 \times 20$ for time-independent problems and $2 \times 20$ for time-dependent problems. \textcolor{black}{In Figure \ref{fig11_ex1}, we display the error versus epoch plots for Example \ref{ch2exmp1}, corresponding to training point values of $\eta_{o} = 1k, 5k$, and $10k$. It is clear that increasing the training data from $1k$ to $10k$ leads to a significant improvement in accuracy. Similarly, Figure \ref{fig111_ex1} shows the loss versus epoch plots for the same example, where the loss represents the sum of the operator loss and the boundary condition loss. These plots correspond to the final iteration of our algorithm, where $(\varepsilon_1, \varepsilon_2) = (\varepsilon_1^{end}, \varepsilon_2^{end}) = (10^{-3}, 10^{-4})$. The plots demonstrate that as the number of training points increases, the loss value decreases substantially. A similar trend was observed in the other examples as well. Figures \ref{err_epoch_ex2} and \ref{loss_ex2} illustrate the results for Example \ref{ch2exmp2}, while Figures \ref{err_epoch_ex2D} and \ref{loss_ex2D} present the corresponding results for Example \ref{ex2D}. A consistent pattern is also evident in Figures \ref{err_epoch_ex2D_time} and \ref{loss_ex2D_time} for Example \ref{ex2D_time}. It was also noted that as the training data dropped below $1k$, the accuracy began to deteriorate. }

\begin{table}[!h]
\centering
\caption{Error comparison using nonlinear update with different widths and depths for Example \ref{ch2exmp1}.}
\label{non_linear_1}
\begin{tabular}{|c|c|c|c|c|}
\hline \text { Depth $\times$ Width } & $(\varepsilon_1,\varepsilon_2)$ & $\mathscr{E}_2$ & CPU time \\
\hline & $(10^{-2}, 10^{-3})$  & 7.38608e-03 & 1379.94 sec \\
$4 \times 20$  & & & \\
 & $(10^{-3}, 10^{-4})$ & 1.21742e-02 & 2137.08 sec\\
\hline & $(10^{-2}, 10^{-3})$  & 7.38585e-03 &  743.09 sec \\
$4 \times 40$  & & & \\
 & $(10^{-3}, 10^{-4})$ & 1.21742e-02 & 1909.18 sec\\
\hline
 & $(10^{-2}, 10^{-3})$  & 7.38600e-03 & 1471.01 sec\\
 $8 \times 20$&  &  &\\
& $(10^{-3}, 10^{-4})$  &  1.21756e-02 & 3138.97 sec\\
\hline
 & $(10^{-2}, 10^{-3})$  & 7.38585e-03 & 1195.67 sec\\
 $8 \times 40$&  &  &\\
& $(10^{-3}, 10^{-4})$  & 1.21736e-02  & 3286.40 sec\\
\hline
\end{tabular}
\end{table}

\begin{table}[!h]
\centering
\caption{Error comparison using nonlinear update with different widths and depths for Example \ref{ch2exmp2}.}
\label{non_linear_2}
\begin{tabular}{|c|c|c|c|c|}
\hline \text { Depth $\times$ Width } & $(\varepsilon_1,\varepsilon_2)$ & $\mathscr{E}_2$ & CPU time \\
\hline & $(10^{-2}, 10^{-3})$  & 3.67415e-03 &  997.53 sec \\
$4 \times 20$  & & & \\
 & $(10^{-3}, 10^{-4})$ & 8.94421e-03 & 2055.41 sec\\
\hline & $(10^{-2}, 10^{-3})$  & 3.67407e-03 & 931.62 sec\\
$4 \times 40$  & & & \\
 & $(10^{-3}, 10^{-4})$ & 8.94421e-03 & 1958.53 sec\\
\hline
 & $(10^{-2}, 10^{-3})$  & 3.67408e-03 & 1291.48 sec\\
 $8 \times 20$&  &  &\\
& $(10^{-3}, 10^{-4})$  &  8.94382e-03 & 3112.42 sec\\
\hline & $(10^{-2}, 10^{-3})$  & 3.67414e-03 & 1571.71 sec  \\
$8 \times 40$  & & & \\
 & $(10^{-3}, 10^{-4})$ & 8.94369e-03 & 3566.46 sec\\
 \hline
\end{tabular}
\end{table}

\begin{table}[!h]
\centering
\caption{Error comparison using linear update with different widths and depths for Example \ref{ch2exmp1}.}
\label{linear_1}
\begin{tabular}{|c|c|c|c|c|}
\hline \text { Depth $\times$ Width } & $(\varepsilon_1,\varepsilon_2)$ & $\mathscr{E}_2$ & CPU time \\
\hline & $(10^{-2}, 10^{-3})$  & 6.67240e-06 &  499.03 sec \\
$4 \times 20$  & & & \\
 & $(10^{-3}, 10^{-4})$ & 5.08541e-04 & 673.43 sec\\
\hline & $(10^{-2}, 10^{-3})$  & 1.42943e-05 &  742.54 sec \\
$4 \times 40$  & & & \\
 & $(10^{-3}, 10^{-4})$ & 2.47156e-04 & 929.41 sec\\
\hline
 & $(10^{-2}, 10^{-3})$  & 5.40406e-06 & 332.67 sec\\
 $8 \times 20$&  &  &\\
& $(10^{-3}, 10^{-4})$  &  4.68982e-04 & 966.12 sec\\
\hline
 & $(10^{-2}, 10^{-3})$  & 2.25027e-05 & 1115.13 sec\\
 $8 \times 40$&  &  &\\
& $(10^{-3}, 10^{-4})$  & 1.404431e-04  & 1615.46 sec\\
\hline
\end{tabular}
\end{table}

\begin{table}[!h]
\centering
\caption{Error comparison using linear update with different widths and depths for Example \ref{ch2exmp2}.}
\label{linear_2}
\begin{tabular}{|c|c|c|c|c|}
\hline \text { Depth $\times$ Width } & $(\varepsilon_1,\varepsilon_2)$ & $\mathscr{E}_2$ & CPU time \\
\hline & $(10^{-2}, 10^{-3})$  & 3.14247e-06 & 365.56 sec \\
$4 \times 20$  & & & \\
 & $(10^{-3}, 10^{-4})$ & 2.13913e-03 & 550.16 sec\\
\hline & $(10^{-2}, 10^{-3})$  & 8.42577e-06 & 346.05 sec \\
$4 \times 40$  & & & \\
 & $(10^{-3}, 10^{-4})$ & 7.52966e-04 & 797.75 sec\\
 \hline
 & $(10^{-2}, 10^{-3})$  & 4.24980e-06 & 482.25 sec\\
 $8 \times 20$&  &  &\\
& $(10^{-3}, 10^{-4})$  &  8.31980e-05 & 961.15 sec\\
\hline
& $(10^{-2}, 10^{-3})$  & 3.18480e-06 & 615.23 sec\\
 $8 \times 40$&  &  &\\
& $(10^{-3}, 10^{-4})$  & 3.87194e-05  & 1405.22 sec\\
\hline
\end{tabular}
\end{table}

\begin{table}[!h]
\centering
\caption{Error comparison using linear update with different widths and depths for Example \ref{ch2exmp3}.}
\label{linear_3}
\begin{tabular}{|c|c|c|c|c|}
\hline \text { Depth $\times$ Width } & $(\varepsilon_1,\varepsilon_2)$ & $\mathscr{E}_2$ & CPU time \\
\hline & $(10^{-2}, 10^{-3})$  & 5.19825e-01 & 442.16 sec \\
$2 \times 20$  & & & \\
 & $(10^{-3}, 10^{-4})$ & 6.28389e-01 & 441.14 sec\\
\hline & $(10^{-2}, 10^{-3})$  & 5.39786e-01 & 648.12 sec \\
$4 \times 20$  & & & \\
 & $(10^{-3}, 10^{-4})$ & 6.60805e-01 & 638.12 sec\\
\hline & $(10^{-2}, 10^{-3})$  & 5.39773e-01 & 807.41 sec \\
$4 \times 40$  & & & \\
 & $(10^{-3}, 10^{-4})$ & 6.51049e-01 & 813.42 sec \\
 \hline
 & $(10^{-2}, 10^{-3})$  & 5.39772e-01 & 1062.02 sec\\
 $8 \times 20$&  &  &\\
& $(10^{-3}, 10^{-4})$  & 6.52699e-01  & 1063.48 sec\\
\hline
& $(10^{-2}, 10^{-3})$  & 5.39773e-01 & 1391.96 sec\\
 $8 \times 40$&  &  &\\
& $(10^{-3}, 10^{-4})$  & 6.28496e-01  & 1407.58 sec\\
\hline
\end{tabular}
\end{table}

\begin{table}[!h]
\centering
\caption{Error comparison using linear update with different widths and depths for Example \ref{ex2D}.}
\label{ex2D_table}
\begin{tabular}{|c|c|c|c|c|}
\hline \text { Depth $\times$ Width } & $(\varepsilon_1,\varepsilon_2)$ & $\mathscr{E}_2$ & CPU time \\
\hline & $(10^{-2}, 10^{-3})$  & 6.87200e-04 & 2756.71 sec \\
$4 \times 20$  & & & \\
 & $(10^{-3}, 10^{-4})$ & 7.37420e-03 & 3608.83 sec\\
\hline & $(10^{-2}, 10^{-3})$  & 4.36345e-04 & 5491.24 sec \\
$4 \times 40$  & & & \\
 & $(10^{-3}, 10^{-4})$ & 6.08680e-03 & 5703.67 sec\\
\hline
 & $(10^{-2}, 10^{-3})$  & 1.94381e-04 & 4363.29 sec\\
 $8 \times 20$&  &  &\\
& $(10^{-3}, 10^{-4})$  & 1.59990e-03 & 6173.14 sec \\
\hline
& $(10^{-2}, 10^{-3})$  & 9.16608e-05 & 7239.75 sec\\
 $8 \times 40$&  &  &\\
& $(10^{-3}, 10^{-4})$  & 1.52576e-03 & 8390.88 sec\\
\hline
\end{tabular}
\end{table}

\subsection{Comparison with some standard numerical methods}
We will compare the results obtained by our proposed algorithm with standard numerical methods, such as the finite difference and finite element methods. In the following tables, the number of mesh points or elements is denoted by $\mathcal{N}$.  Table \ref{table_comp1} provides a comparison between the results of the proposed method and the Upwind Finite Difference Method (Upwind FDM) for Example \ref{ch2exmp1}, confirming that PA-PINNs achieve higher accuracy than Upwind FDM. Additionally, Table \ref{table_comp2} presents a comparison of our results with both Upwind FDM and the standard Finite Element Method (FEM) for Example \ref{ch2exmp2}, clearly demonstrating that PA-PINNs outperform both methods in terms of accuracy. In Table \ref{table_comp3}, a comparison for Example \ref{ex2D} with the standard FEM is provided, showing that PA-PINNs offer better results. Same trends were evident for other examples as well. It is well known that numerical methods like finite difference and finite element methods face several challenges when applied to two-parameter SPPs. These methods often require intricate mesh refinement or specialized techniques to capture boundary layers, like implementation over layer adapted meshes \cite{2010_Linss}(eg., Shishkin/Bakhvalov/Vulanovic meshes, etc) for capturing the boundary layers effectively, but that again depends on prior knowledge of the boundary layer's width and location, making the implementation cumbersome. Also, two-parameter SPPs exhibit varying solution behaviors depending on the ratio of the two parameters as discussed in Section \ref{model_problem_two_para} , further complicating the implementation of these numerical methods. In contrast, our proposed PA-PINNs approach adapts naturally to different parameter regimes without needing any modification to the network architecture, simplifying implementation while offering better accuracy, as demonstrated by our error comparisons with finite difference and finite element methods.

\begin{table}[!h]
\centering
\caption{\textcolor{black}{Error comparison for Example \ref{ch2exmp1}.}}
\label{table_comp1}
\begin{tabular}{|c|c|c|c|c|}
\hline \text {\textcolor{black}{Method}} & $(\varepsilon_1,\varepsilon_2)$ & $\mathscr{E}_{\infty}$ \\
\hline & $(10^{-2}, 10^{-3})$  & 1.18356e-05 \\
\textcolor{black}{Upwind FDM} & &  \\
\textcolor{black}{(over Shishkin mesh)} & $(10^{-3}, 10^{-4})$ & 7.68553e-03 \\
\textcolor{black}{with $\mathcal{N}$ = 1024} &&\\
& $(10^{-4}, 10^{-5})$ & 6.62190e-02\\
\hline
& $(10^{-2}, 10^{-3})$  & 1.26608e-05\\
\textcolor{black}{Proposed method} &  &  \\
\textcolor{black}{with depth $\times$ width $= 8 \times 20$}& $(10^{-3}, 10^{-4})$  & 9.75883e-04\\
\textcolor{black}{and $\eta_{o} = 1000$} &&\\
& $(10^{-4}, 10^{-5})$ & 4.50602e-03\\
\hline
\end{tabular}
\end{table}

\begin{table}[!h]
\centering
\caption{\textcolor{black}{Error comparison for Example \ref{ch2exmp2}.}}
\label{table_comp2}
\begin{tabular}{|c|c|c|c|c|}
\hline \text {\textcolor{black}{Method}} & $(\varepsilon_1,\varepsilon_2)$ & $\mathscr{E}_{\infty}$ \\
\hline & $(10^{-1}, 10^{-2})$ & 1.13160e-02\\
\textcolor{black}{Upwind FDM}&&\\
\textcolor{black}{(over Shishkin mesh)}& $(10^{-2}, 10^{-3})$  & 7.27137e-03 \\
\textcolor{black}{with $\mathcal{N}$ = 1024} & &  \\
 & $(10^{-3}, 10^{-4})$ & 7.72630e-03 \\
\hline & $(10^{-1}, 10^{-2})$ & 1.13324e-02\\
\textcolor{black}{Standard FEM} &&\\
\textcolor{black}{(over Shishkin mesh)}& $(10^{-2}, 10^{-3})$  & 7.24019e-03\\
\textcolor{black}{with $\mathcal{N}$ = 1024} &  &  \\
& $(10^{-3}, 10^{-4})$  & 2.81312e-03\\
\hline & $(10^{-1}, 10^{-2})$ & 5.08473e-03\\
\textcolor{black}{Proposed method} &&\\
\textcolor{black}{with depth $\times$ width $= 8 \times 20$}& $(10^{-2}, 10^{-3})$  & 9.27661e-06\\
\textcolor{black}{and $\eta_{o} = 1000$}  &  &  \\
& $(10^{-3}, 10^{-4})$  &  6.79974e-04\\
\hline
\end{tabular}
\end{table}


\begin{table}[!h]
\centering
\caption{\textcolor{black}{Error comparison for Example \ref{ex2D}.}}
\label{table_comp3}
\begin{tabular}{|c|c|c|c|c|}
\hline \text {\textcolor{black}{Method}} & $(\varepsilon_1,\varepsilon_2)$ & $\mathscr{E}_{\infty}$ \\
\hline & $(10^{-2}, 10^{-3})$  & 6.68075e-04 \\
\textcolor{black}{Standard FEM} & &  \\
\textcolor{black}{(over Shishkin mesh)} & $(10^{-3}, 10^{-4})$ & 6.87868e-03 \\
\textcolor{black}{with $\mathcal{N} = 64 \times 64 \times 2$} &&\\
\hline
& $(10^{-2}, 10^{-3})$  & 1.94381e-04 \\
\textcolor{black}{Proposed method} &  &  \\
\textcolor{black}{with depth $\times$ width $= 8 \times 20$}& $(10^{-3}, 10^{-4})$  & 1.59990e-03\\
\textcolor{black}{and $\eta_{o} = 10k$} &&\\
\hline
\end{tabular}
\end{table}

\section{Conclusion}\label{conclusion}
 In this article, PA-PINNs strategy has been implemented for two-parameter elliptic and parabolic SPPs in one- and two-dimensions. Such problems appear in chemical flow reactor theory, geophysics and many other fields of engineering and sciences. By taking a larger value of the perturbation parameter, the smooth part is approximated first, and then gradually reducing the perturbation parameter approximates the layer parts. This strategy enables PINNs to efficiently predict the singular solution. Numerical experiments show that PA-PINNs is effective in solving two-parameter SPPs without a priori knowledge of boundary layer location and width, produces stable and convergent approximations of the solutions. As the layer structure for such problems depends on the ratio of $\varepsilon_2^2/\varepsilon_1$, standard numerical techniques require different mesh discretization and analysis for different cases. However, PA-PINNs can tackle these cases as a whole and hence no change in the architecture. Our tests show that PA-PINNs outperformed standard FDM and FEM in terms of accuracy. By our previous studies on two-paramter SPPs using finite basis PINNs (FB-PINNs) \cite{raina2023comparative}, we can conclude that PA-PINNs is superior than FB-PINNs in terms of accuracy and convergence for problems involving thin boundary layers.

\goodbreak\noindent
\section*{Declarations}

\subsection*{Acknowledgment}
The second author would like to thank Indian Institute of Technology Guwahati for supporting him financially by giving the fellowship for his research work.

\subsection*{Funding}
Not available.

\subsection*{Data Availability}
The code is made available upon reasonable request.

\subsection*{Conflict of interest}
The authors declare that they have no conflict of interest.

\subsection*{Ethical approval}
This article does not contain any studies with human participants or animals performed by any of the authors.




\begin{thebibliography}{}

  
\bibitem{arzani2023theory}
{\sc A.~Arzani, K.~W. Cassel, and R.~M. D'Souza}, {\em Theory-guided
  physics-informed neural networks for boundary layer problems with singular
  perturbation}, Journal of Computational Physics, 473 (2023), p.~111768.
  
  
\bibitem{bajaj2023recipes}
{\sc C.~Bajaj, L.~McLennan, T.~Andeen and A.~Roy}, {\em Recipes for when physics fails: recovering robust learning of physics informed neural networks}, Machine learning: science and technology, 4 (2023), p.~015013. 
  
  
\bibitem{barman2023alternating}
{\sc M.~Barman, S.~Natesan, and A.~Sendur}, {\em Alternating direction implicit method
  for singularly perturbed 2{D} parabolic convection--diffusion--reaction problem
  with two small parameters}, International Journal of Computer Mathematics, 100 (2023), p.~253--282. 
  
\bibitem{barman2024parameter}
{\sc M.~Barman, S.~Natesan, and A.~Sendur}, {\em A parameter-uniform hybrid method for singularly perturbed parabolic 2D convection-diffusion-reaction problems}, Applied Numerical Mathematics, 207 (2025), p.~111--135.   
  
  
\bibitem{brdar2016singularly}
{\sc M.~Brdar, H.~Zarin}, {\em A singularly perturbed problem with two parameters on a Bakhvalov-type mesh }, Journal of Computational and Applied Mathematics, 292 (2016), p.~307--319.    


\bibitem{cao2023physics}
{\sc F.~Cao, F.~Gao, X.~Guo, and D.~Yuan}, {\em Physics-informed neural
  networks with parameter asymptotic strategy for learning singularly perturbed convection-dominated problem},
  Computers \& Mathematics with Applications, 150 (2023), pp.~229--242.

\bibitem{chen1974asymptotic}
{\sc J.~Chen and RE.~O’Malley Jr}, {\em On the asymptotic solution of a two-parameter boundary value problem of chemical reactor theory}, SIAM Journal on Applied Mathematics, 26 (1974), pp.~717--729.

\bibitem{cheng2021local}
{\sc Y.~Cheng}, {\em On the local discontinuous Galerkin method for singularly perturbed problem with two parameters}, Journal of Computational and Applied Mathematics, 392 (2021), pp.~113485.

\bibitem{avijit2022convergence}
{\sc A.~Das and S.~Natesan}, {\em Convergence analysis of a fully-discrete FEM for singularly perturbed two-parameter parabolic PDE}, Mathematics and Computers in Simulation, 197 (2022), pp.~185--206.


\bibitem{diprima1968asymptotic}
{\sc R.C.~DiPrima}, {\em Asymptotic methods for an infinitely long slider squeeze-film bearing}, Journal of Lubrication Technology, 90 (1968), pp.~173--183.

\bibitem{farrel_book_2000}
{\sc P.~Farrell, A.~Hegarty, J.~Miller, E.~O'Riordan, G.~Shishkin}, Robust
	Computational Techniques for Boundary Layers, Chapman {\&} Hall/CRC Press,
	Boco Raton, 2000.

\bibitem{franz2014c0}
{\sc S.~Franz, H.-G. Roos, and A.~Wachtel}, {\em A ${C^{0}}$ interior penalty
  method for a singularly-perturbed fourth-order elliptic problem on a
  layer-adapted mesh}, Numerical Methods for Partial Differential Equations, 30
  (2014), pp.~838--861.
  
\bibitem{gao2021phygeonet}
{\sc H.~Gao, L.~Sun, and J.X.~Wang}, {\em PhyGeoNet: Physics-informed geometry-adaptive convolutional neural networks for solving parameterized steady-state PDEs on irregular domain}, Journal of Computational Physics, 428 (2021), pp.~110079.  

\bibitem{hemker2000varepsilon}
{\sc P.W.~Hemker, G.I.~Shishkin, and L.P.~Shishkina}, {\em $\varepsilon$-uniform schemes with high-order time-accuracy for parabolic singular perturbation problems}, IMA Journal of Numerical Analysis, 20 (2000), pp.~99--121.  
  
\bibitem{hornik1991approximation}
{\sc K.~Hornik}, {\em  Approximation capabilities of multilayer feedforward networks}, Neural networks, 4 (1991), pp.~251--257.    
 

\bibitem{jagtap2021extended}
{\sc A.~D. Jagtap and G.~E. Karniadakis}, {\em Extended physics-informed neural
  networks (XPINNS): A generalized space-time domain decomposition based deep
  learning framework for nonlinear partial differential equations.}, in AAAI
  spring symposium: MLPS, vol.~10, 2021.
  
  
\bibitem{jagtap2020adaptive}
{\sc A.~D. Jagtap, K.~Kawaguchi, and G.~E. Karniadakis}, {\em Adaptive activation functions accelerate convergence in deep and physics-informed neural networks},Journal of Computational Physics, 404 (2020), p.~109136.  
  
  
\bibitem{jagtap2020conservative}
{\sc A.~D. Jagtap, E.~Kharazmi, and G.~E. Karniadakis}, {\em Conservative
  physics-informed neural networks on discrete domains for conservation laws:
  Applications to forward and inverse problems}, Computer Methods in Applied
  Mechanics and Engineering, 365 (2020), p.~113028.
  
\bibitem{karniadakis2021physics}
{\sc G.~E. Karniadakis, I.G.~Kevrekidis, L.~Lu, P.~Perdikaris, S.~Wang and L.~Yang}, {\em Physics-informed machine learning}, Nature Reviews Physics, 3 (2021), p.~422--440.   
  
  
\bibitem{kharazmi2021hp}
{\sc E.~Kharazmi, Z.~Zhang, and G.~E. Karniadakis}, {\em $hp$-VPINNS: Variational
  physics-informed neural networks with domain decomposition}, Computer Methods
  in Applied Mechanics and Engineering, 374 (2021), p.~113547.  
  
\bibitem{2010_Linss}
{\sc T.~Lin{$\ss$}}, Layer-adapted meshes for reaction-convection-diffusion problems, Springer, 2009.    
  
\bibitem{linss2004analysis}
{\sc T.~Lin{$\ss$} and H.G.~Roos}, {\em Analysis of a finite-difference scheme for a singularly perturbed problem with two small parameters}, Journal of Mathematical Analysis and Applications, 289 (2004), p.~355--366.  


\bibitem{majumdar2017alternating}
{\sc A.~Majumdar and S.~Natesan}, {\em Alternating direction numerical scheme for singularly perturbed 2D degenerate parabolic convection-diffusion problems}, Applied Mathematics and Computation, 313 (2017), p.~453--473.  
  

\bibitem{mcclenny2023self}
{\sc L.~D. McClenny and U.~M. Braga-Neto}, {\em Self-adaptive physics-informed neural networks}, Journal of Computational Physics, 474 (2023), p.~111722.  
  
\bibitem{miller_book_1996}
{\sc J.~Miller, E.~O'Riordan, G.~Shishkin}, Fitted Numerical Methods for Singular
	Perturbation Problems, World Scientific, Singapore, 1996.
	
\bibitem{mohapatra2010parameter}
{\sc J.~Mohapatra and S.~Natesan}, {\em Parameter-uniform numerical method for global solution and global normalized flux of singularly perturbed boundary value problems using grid equidistribution}, Computers \& Mathematics with Applications, 60 (2010), p.~1924--1939.	
	
	
\bibitem{moseley2023finite}
{\sc B.~Moseley, A.~Markham, and T.~Nissen-Meyer}, {\em Finite basis
  physics-informed neural networks (FBPINNS): a scalable domain decomposition
  approach for solving differential equations}, Advances in Computational
  Mathematics, 49 (2023), p.~62.
  
\bibitem{o1967two}
{\sc R.E~O'Malley}, {\em Two-parameter singular perturbation problems for second-order equations}, Journal of Mathematics and Mechanics, 16 (1967), p.~1143--1164.  
   

\bibitem{malley_book_1974}
{\sc RE.~O'Malley}, Introduction to singular perturbations, Academic Press, New york (1974).

\bibitem{malley_book_1991}
{\sc RE.~O'Malley Jr.}, Singular Perturbation Methods for Ordinary Differential Equations, Springer, New York (1991).

\bibitem{o2011parameter}
{\sc E.~O'Riordan, and M.L.~Pickett}, {\em A parameter-uniform numerical method for a singularly perturbed two parameter elliptic problem}, Advances in Computational Mathematics, 35 (2011), p.~57--82.

\bibitem{o2006numerical}
{\sc E.~O'Riordan, M.L.~Pickett, and G.I~Shishkin}, {\em Numerical methods for singularly perturbed elliptic problems containing two perturbation parameters}, Mathematical Modelling and Analysis, Taylor \& Francis, 11 (2006), p.~199--212.


 
\bibitem{pang2019fpinns}
{\sc G.~Pang, L.~Lu, and G.~E. Karniadakis}, {\em fPINNS: Fractional
  physics-informed neural networks}, SIAM Journal on Scientific Computing, 41
  (2019), pp.~A2603--A2626.


\bibitem{raissi2019physics}
{\sc M.~Raissi, P.~Perdikaris, and G.~E. Karniadakis}, {\em Physics-informed
  neural networks: A deep learning framework for solving forward and inverse
  problems involving nonlinear partial differential equations}, Journal of
  Computational physics, 378 (2019), pp.~686--707.

\bibitem{raina2023comparative}
{\sc A.~Raina, B.~Satyadev, and S.~Natesan}, {\em A comparative study on solving two-parameter singular perturbation problems using recent physics inspired neural networks, Submitted for publication.}

\bibitem{raina2024fourth}
{\sc A.~Raina, and S.~Natesan}, {\em A weak Galerkin finite element method for fourth-order parabolic singularly perturbed problems on layer adapted Shishkin mesh}, Applied Numerical Mathematics, 207 (2025), pp.~520--533.

\bibitem{raina2024wg}
{\sc A.~Raina, S.~Natesan, and \c{S}.~Toprakseven}, {\em Anisotropic error analysis of weak Galerkin finite element method for singularly perturbed biharmonic problems}, Mathematics and Computers in Simulation, 229 (2025), pp.~203--221.

\bibitem{raina2025novel}
{\sc A.~Raina, S.~Natesan, and \c{S}.~Toprakseven}, {\em A Novel ADI-Weak Galerkin Method for Singularly Perturbed Two-Parameter 2D Parabolic PDEs}, Numerical Methods for Partial Differential Equations, 41(1) (2025), pp.~e23169.




\bibitem{roos_book_2008}
{\sc H.-G. Roos, M.~Stynes, L.~Tobiska}, Robust Numerical Methods for Singularly
	Perturbed Differential Equations: Convection-Diffusion-Reaction and Flow
	Problems, Vol.~24, Springer Science \& Business Media, 2008.
	
\bibitem{shukla2021parallel}
{\sc K.~Shukla, A.D.~Jagtap and G.~E. Karniadakis}, {\em Parallel physics-informed neural networks via domain decomposition}, Journal of Computational Physics, 447
  (2021), pp.~110683.	
	
\bibitem{singh2020study}
{\sc G.~Singh and S.~Natesan}, {\em Study of the NIPG method for two--parameter singular perturbation problems on several layer adapted grids}, Journal of Applied Mathematics and Computing, 63 (2020), pp.~683--705.	
	

\bibitem{singh2020parameter}
{\sc M.~K. Singh and S.~Natesan}, {\em A parameter-uniform hybrid finite
  difference scheme for singularly perturbed system of parabolic
  convection-diffusion problems}, International Journal of Computer
  Mathematics, 97 (2020), pp.~875--905.


\bibitem{teofanov2007elliptic}
{\sc L.J.~Teofanov and H.G.~Roos}, {\em An elliptic singularly perturbed problem with two parameters I: Solution decomposition}, Journal of Computational and Applied Mathematics, 206 (2007), pp.~1082--1097.



\bibitem{toprakseven2023weak}
{\sc {\c{S}}.~Toprakseven, A.~Kaushik, and M.~Sharma}, {\em A weak Galerkin finite
  element method for singularly perturbed problems with two small parameters on
  bakhvalov-type meshes}, Numerical Algorithms,  (2023), pp.~1--25.


\bibitem{yu2022gradient}
{\sc J.~Yu, L.~Lu, X.~Meng, and G.~E. Karniadakis}, {\em Gradient-enhanced
  physics-informed neural networks for forward and inverse PDE problems},
  Computer Methods in Applied Mechanics and Engineering, 393 (2022), p.~114823.

\bibitem{zhang2020learning}
{\sc D.~Zhang, L.~Guo, and G.~E. Karniadakis}, {\em Learning in modal space:
  Solving time-dependent stochastic PDEs using physics-informed neural
  networks}, SIAM Journal on Scientific Computing, 42 (2020), pp.~A639--A665.

\bibitem{zhang2021high}
{\sc J.~Zhang and Y.~Lv}, {\em High-order finite element method on a
  {B}akhvalov-type mesh for a singularly perturbed convection--diffusion problem with two parameters}, Applied Mathematics and Computation, 397 (2021),
  p.~125953.



\end{thebibliography}
\end{document}